\documentclass[10pt]{article}
\usepackage{amsfonts}
\usepackage{mathrsfs}
\usepackage{amsmath,amsthm,amsfonts,amssymb,epsfig}

\newtheorem{thm}{Theorem}[section]

\newtheorem{lem}[thm]{Lemma}

\theoremstyle{definition}

\theoremstyle{remark}

\numberwithin{equation}{section}

\newcommand\be{\begin{equation}}
\newcommand\ee{\end{equation}}
\newcommand\bes{\begin{eqnarray}}
\newcommand\ees{\end{eqnarray}}
\newcommand\bess{\begin{eqnarray*}}
\newcommand\eess{\end{eqnarray*}}

\begin{document}
\title{Carleman estimates for space semi-discrete approximations of one-dimensional stochastic parabolic equation and its applications}

\author{ Bin Wu$^{1,2}$\footnote{Corresponding author. email: binwu@nuist.edu.cn}\quad Ying Wang$^1$\quad Zewen Wang$^3$
   \\
$ ^1$School of Mathematics and Statistics\\ Nanjing University of
Information Science
and Technology\\
Nanjing 210044, China\\
$ ^2$Center for Applied Mathematics of Jiangsu Province \\
Nanjing University of Information Science and
Technology \\
Nanjing 210044, China\\
$ ^3$ Department of Basic Courses\\
Guangzhou Maritime University\\
Guangdong 510725, China
}

\maketitle

\begin{abstract}
In this paper, we study discrete Carleman estimates for space semi-discrete approximations of one-dimensional stochastic parabolic equation. As applications of these discrete Carleman estimates, we apply them to study two inverse problems for the spatial semi-discrete stochastic parabolic equations, including a discrete inverse random source problem and a discrete Cauchy problem. We firstly establish two Carleman estimates for a one-dimensional semi-discrete stochastic parabolic equation, one for homogeneous boundary and the other for non-homogeneous boundary. Then we apply these two estimates separately to derive two stability results. The first one is the Lipschitz stability for the discrete inverse random source problem. The second one is the H\"{o}lder stability for the discrete Cauchy problem.

\vskip 0.3cm


\end{abstract}

\section{Introduction}
\setcounter{equation}{0}

Let $(\Omega, \mathcal{F}, \{\mathcal{F}_t\}_{t\geq0}, \mathbb P)$ be a complete filtered probability space on which a one-dimensional standard Brownian motion $\{B(t)\}_{t\geq0}$ is defined such that $\{\mathcal F_t\}_{t\geq 0}$ is the natural filtration generated by $B(\cdot)$, augmented by all the $\mathbb P$-null sets in $\mathcal{F}$.  Let $T>0$, $G=(0,L)$ and $Q=G\times(0,T)$. Consider the following one-dimensional stochastic parabolic equation
\begin{align}\label{1.1}
	\left\{
	\begin{array}{ll}
		{\rm d}y-y_{xx}{\rm d}t=(a y+ b y_x+f){\rm d}t+(cy+g){\rm d} B(t), &(x,t)\in Q,\\
        y(0,t)=y(L,t)=0, 	&t\in (0,T), \\
	    y(x,0)=y^0(x),		    & x\in G
	\end{array}
	\right.
\end{align}
with suitable coefficients $a,b$ and $c$.  Physically, $f$ and $g$ are source terms, $g$ stands for the intensity of a random force of the white noise
type.

The main objective of this paper is to establish discrete Carleman estimates for the finite-difference space discretization of stochastic parabolic equation (\ref{1.1}), and to give its applications in two kinds of typical inverse problems, namely the inverse random source problem and the Cauchy problem for spatial semi-discrete stochastic parabolic equation. To do this, let us consider $N\in \mathbb{N}^*$, a step $h=\frac{L}{N+1}$, and
an equidistant mesh of the interval $(0, L), 0=x_0<x_1<\cdots<x_N<x_{N+1}=L$, with
$x_i=ih, 0\leq i\leq N + 1$. The finite-difference approximation of the space
derivatives leads to the following semi-discretization of (1.1):
\begin{equation}\label{1.2}
	\left\{
	\begin{array}{l}
		{\rm d}y_i(t)-\frac{y_{i+1}(t)-2y_i(t)+y_{i-1}(t)}{h^2}{\rm d}t=\left(a_iy_i(t)+b_i\frac{y_{i+1}(t)-y_i(t)}{h}+f_i\right){\rm d}t\\
\hspace{5cm}+(c_iy_i(t)+g_i){\rm d}B(t), \hspace{0.6cm} 1\leq i\leq N,\ t>0,\\
		y_{0}(t)=y_{N+1}(t)=0, \hspace{5.7cm} t>0,\\
		y_{i}(0)=y^0_i,  \hspace{7.15cm}   1\leq i\leq N,
	\end{array}
	\right.
\end{equation}
 where $y_i(t)$ stands for $y(x_i,t)$, analogous definitions are for $a_i, b_i, f_i, g_i$ and the other functions in the sequel. 
 
 Throughout this paper, we denote the discrete domains and discrete boundaries by the following notations
\begin{equation*}
	\begin{aligned}
		&G_{h}=\left\{x_1, x_2,\cdots, x_N\right\}, \quad&&Q_{h}=G_{h}\times(0,T),\\ 
		&\overline{G}_{h}=\left\{x_0, x_1, x_2, \cdots,  x_N, x_{N+1} \right\}, \quad&& \overline{Q}_{h}= \overline{G}_{h}\times(0,T),\\ 
		&\partial G_{h} =\left\{x_0, x_{N+1}\right\}, \quad&&\Sigma_{h}=\partial G_{h}\times(0,T),\\ 
		&G_{h}^{-}=\left\{x_0, x_1, x_2, \cdots,  x_N\right\}, \quad&&Q_{h}^{-}=G_{h}^{-}\times(0,T).
	\end{aligned}
\end{equation*}
We denote by $\mathbb{R}^{\mathfrak M}$ and $\mathbb{R}^ {\overline{\mathfrak M}}$ the set of discrete functions defined on $G_h$ and $\overline{G}_h$, respectively. Furthermore, for $u_h=(u_0,u_1,\cdots,u_{N+1})^{\rm T}\in \mathbb{R}^ {\overline{\mathfrak M}}$, we define the averaging operators and difference operators  as follows
\begin{equation*}
\begin{aligned}
	&(\mathbf{m}_h^{+}u_h)_{i}=\frac{u_{i+1}+u_{i}}{2},\quad (\mathbf{m}_{h}^{-}u_h)_{i}=\frac{u_{i}+u_{i-1}}{2},\quad (\mathbf{m}_h u_h)_{i}=\frac{u_{i+1}+2u_{i}+u_{i-1}}{4},\\
&(\mathbf{D} _{h}^{+}u_h)_{i}=\frac{u_{i+1}-u_{i}}{h},\quad(\mathbf{D}_{h}^{-}u_h)_{i}=\frac{u_{i}-u_{i-1}}{h}, \quad(\mathbf{D} _{h}u_h)_{i}=\frac{u_{i+1}-u_{i-1}}{2h},\\
&(\Delta _{h}u_h)_{i}=\frac{u_{i+1}-2u_{i}+u_{i-1}}{h^{2}}.
\end{aligned}
\end{equation*}
These notations allow us to express semi-discretization (\ref{1.2}) in a more compact way, which is necessary to formulate our inverse problems.

\vspace{2mm}
 
\noindent{\bf Discrete inverse random source problem.} \ Determine the random source term $g_h\in \mathbb{R}^{\mathfrak M}$ in the following semi-discrete stochastic parabolic equation:
\begin{align}\label{1.3}
	\left\{\begin{array}{ll}
		{\rm d} y_h-\Delta_{h}y_h{\rm d}t=(a_hy_h+b_h{\mathbf D}^+_h y_h){\rm d}t+g_h{\rm d}B(t), &(x_h,t)\in Q_{h},\\
		y_h=0, &(x_h,t)\in \Sigma_{h}, \\
        y_h(0)=y_h^0, &x_h\in G_{h},
        \end{array}
	\right.
\end{align}
by the observation data
 $$(\mathbf D_{h}^{-}y_h)_{N+1}\quad{\rm and}\quad y_h(T),$$
where $y_h\in \mathbb{R}^{\overline{\mathfrak M}}$ and $a_h, b_h, y_h^0\in \mathbb{R}^{\mathfrak M}$.

\vspace{2mm}

Let $G_0\subset G$ such that $\overline {G_0}\subset G\cup \{x_{N+1}\}$ and $\partial G_0\cap \partial G=\{x_{N+1}\}$. We set $G_{0,h}=G_0\cap G_h$. The discrete Cauchy problem is described as follows.

\vspace{2mm}

\noindent{\bf Discrete Cauchy problem.} For any $\epsilon>0$, determine the solution $y_h$ in $G_{0,h}\times(\epsilon, T-\epsilon)$ of the following semi-discrete stochastic parabolic equation:
\begin{equation}\label{1.4}
		{\rm d}y_{h}-\Delta _{h}y_{h}{\rm d}t=(a_hy_{h}+b_h{\mathbf D}_{h}^{+}y_{h}){\rm d}t+c_hy_{h}{\rm d}B(t),\quad  (x_h,t)\in Q_{h},
\end{equation}
by the discrete lateral boundary data
\begin{equation}
(y_{h})_{N+1}=\xi(t) \quad {\rm and} \quad \left(\mathbf D_{h}^{-} y_{h}\right)_{N+1}=\eta(t), \quad t\in (0,T).
\end{equation}

\vspace{2mm} 
 
 Carleman estimate is a class of weighted energy estimates related to some differential operator, which can be applied to many aspects, such as inverse problems [\ref{Imanuvilov2005}, \ref{Kli1992}, \ref{Kli2004}, \ref{Ya2009}], control theory [\ref{Fernández2006}, \ref{Isakov2006}, \ref{Rousseau2010}] and so on.  
 There are rich references on discrete Carleman estimates for deterministic partial differential equations, which can be divided into three main categories: space-discrete, time-discrete and full-discrete results. We refer to [\ref{Bau2015}, \ref{Bau2013}, \ref{Boy2010}, \ref{Boyer2010}, \ref{Boy2014}, \ref{Cer2022}, \ref{Lec2021}, \ref{Ngu2012}] for  space-discrete Carleman estimates, [\ref{Boy2020}, \ref{Hernández2021}] for time-discrete Carleman estimates and [\ref{Boy2011}, \ref{Cas2021}] for full-discrete Carleman estimates. These discrete Carleman estimates have been successfully applied to prove the uniform controllability and the stability of inverse problems for various discrete deterministic partial differential equations [\ref{Lec2023}, \ref{Per2022}, \ref{Zhang2022}, \ref{Zhao2023}].  To the best of our knowledge, there is no paper considering discrete Carleman estimates for stochastic
partial differential equations.

Carleman estimate is a powerful tool to study inverse problems related to various stochastic partial differential equations [\ref{Gao2018}, \ref{Liu2019}, \ref{Wu2020}, \ref{Wu2022}, \ref{Yuan2017}]. For inverse source problem related to stochastic parabolic equations, we refer to [\ref{Lu2012}] for the uniqueness of an inverse problem of determining source function $f$ in (\ref{1.1}).   The inverse source problem of determining two kinds of sources $f$ and $g$ simultaneously in (\ref{1.1}) was studied in [\ref{Yuan2021}].  Moreover, we also refer to  [\ref{Bao2010}] or [\ref{Bao2012}] for applications of regularization techniques in the numerical methods for inverse random source problems. The Cauchy problem aims to recover the solution with observed data from the lateral boundary. 
In [\ref{Ya2009}], a conditional stability was proved for the Cauchy problem of deterministic parabolic equations.  Recently, this method is extended to the stochastic case. The conditional stability and convergence rate of the Tikhonov regularization method for the Cauchy problem of stochastic parabolic equations was obtained [\ref{Dou2023}]. It is worth mentioning that in [\ref{Lu2023}] the authors gave a detailed review on inverse problems for stochastic partial differential equations. However, these results were all obtained within a continuous framework. Discrete inverse problems for the stochastic differential equations have not been studied thoroughly yet.

In this paper, we firstly focus on Carleman estimates for discrete stochastic parabolic equation. More precisely, we will prove two Carleman estimates for spatial semi-discretization (\ref{1.2}) of one-dimensional stochastic parabolic equation.  We apply the first Carleman estimate to study the discrete inverse random source problem. Unlike the deterministic counterparts, the solution of a stochastic differential equation is not differentiable with respect to time variable. We have to choose a regular weight function to put the random source term on the left-hand side of discrete Carleman estimate. Since the second parameter $\lambda$ plays an important role in the proof of the uniform stability result with respect to the mesh size, we need to carefully decouple $\lambda$ from the constant $C$ in our Carleman estimate. 
 Secondly, in order to handle the discrete Cauchy problem, we need a slight revision to present the second discrete Carleman estimate with non-homogeneous boundary. Applying this Carleman estimate, we obtain a H\"{o}lder stability for the discrete Cauchy problem. In comparison with the deterministic discrete Carleman estimates, there are additional terms arising from stochastic effects to be considered.

 The rest of this paper is organized as follows. In section 2, we present discrete settings and our main results. In section 3, we show two Carleman estimates for space semi-discrete approximations of one-dimensional stochastic parabolic equation. In next two sections, based on these two Carleman estimates we study
the discrete  inverse random source problem and the discrete Cauchy problem, respectively.

\section{Discrete settings and main results}
\setcounter{equation}{0}

In this section, we introduce fundamental concepts of discrete calculus, including discrete function spaces, some useful discrete identities, and integration by parts for discrete operators, which will be used in the proofs of our main results.  Subsequently, we present our main results in this paper. We first give two discrete Carleman estimates for the semi-discrete stochastic parabolic equation, namely, one for homogeneous boundary and the other for non-homogeneous boundary. Then we address two stability results. The first one pertains to the discrete inverse random source problem, while the second one focuses on the discrete Cauchy problem.

\subsection{Discrete settings}
By analogy with the continuous case, for $u_h=(u_0,u_1,\cdots,u_{N+1})^{\rm T}\in \mathbb{R}^{\overline{\mathfrak M}}$ we define the discrete  integrals:\\
\begin{equation*}
	\begin{aligned}
		&\int_{ G_{h} }u_{h}=h\sum_{x_i\in{G_{h}}}u_{i}=h\sum_{i=1}^{N}u_{i},\quad\int_{ G_{h}^{-} }u_{h}=h\sum_{x_i\in{G_{h}^{-}}}u_{i}=h\sum_{i=0}^{N}u_{i}.
	\end{aligned}
\end{equation*}
For $u_h, v_h\in \mathbb{R}^{\mathfrak M}$, we define the following $L^2$-inner product on $\mathbb{R}^{\mathfrak M}$ 
\begin{align*}
(u_h,v_h)_{L^2(G_h)}=h\sum_{i=1}^{N}u_{i}v_i.
\end{align*}
The associated norm is denoted by  $\|u_h\|_{L^2(G_h)}$.  Analogously, we define the $L^\infty$-norm on $\mathbb{R}^{\mathfrak M}$
\begin{equation*}
	\begin{aligned}
		||u_{h}||^2_{L^{\infty}(G_{h})}=\max_{1\leq i\leq N}|u_i|.
	\end{aligned}
\end{equation*}
Also, we introduce
\begin{align*}
\|u_h\|^2_{H^1(G_h)}=\int_{G^{-}_h} |\mathbf D^{+}_h u_h|^2+\int_{G_h}|u_h|^2
\end{align*}
and
\begin{align*}
\|u_h\|_{H^2(G_h)}=\int_{G_h} |\Delta_h u_h|^2+\int_{G^{-}_h} |\mathbf D^{+}_h u_h|^2+\int_{G_h}|u_h|^2.
\end{align*}
Furthermore, for a discrete Banach space $\mathcal X_h$ defined on ${G_h}$,  we denote by $L^2(0,T;\mathcal X_h)$ the set of discrete functions endowed with the norm
\begin{equation*}
			||u_{h}||^2_{L^{2}(0,T;\mathcal X_h)}=\int_{0}^T\|u_h\|_{\mathcal X_h}^2{\rm d}t.
\end{equation*}

Now we introduce some notations for stochastic analysis on discrete space meshes. For a Banach space $\mathcal Y$, we denote by $L_\mathcal F^2(\Omega;\mathcal Y)$ the space of all progressively measurable stochastic process $\zeta$ such that $\mathbb E(\|\zeta\|^2_{\mathcal Y})<\infty$.  For a discrete Banach space $\mathcal X_h$ on ${G_h}$, we denote by  $L^2_{\mathcal{F}}(0, T; \mathcal X_h)$ the Banach space consisting of all $\mathcal X_h$-valued $\{\mathcal F_t\}_{t\geq0}$-adapted processes $\zeta(\cdot)$ such that $\mathbb{E}(\|\zeta(\cdot)\|^2_{L^2(0,T;\mathcal X_h)}) $ $<\infty$, with the canonical norm; by $L^\infty_\mathcal{F}(0, T; \mathcal X_h)$ the Banach space consisting of all $\mathcal X_h$-valued $\{\mathcal F_t\}_{t\geq 0}$-adapted bounded processes.

 For the space-discrete operators, we also need to provide several preliminary identities and discrete integration by parts formula that will be extensively used in the sequel. We present the results without the proof and refer the readers to [\ref{Boy2010}]  (see also [\ref{Bau2013}]) for a detailed discussion.

\begin{lem} Let $ u$ and $v$ be discrete functions defined on $\overline G_h$. Then for the averaging operators and
difference operators, we have the following identities:
\begin{align}
\label{2.1}& {\mathbf m}_{h}^{+}u=u+\frac{h}{2} {\mathbf D}_{h}^{+}u,\   {\mathbf m}_{h}u=u+\frac{h^{2}}{4}\Delta _{h}u,\  \mathbf D_h u=\mathbf m_h^-\mathbf D_h^+u,\  \Delta_{h}u=\mathbf D_{h}^{+}\mathbf D_{h}^{-}u,\\
&\mathbf m_{h}^{+}(uv)=\mathbf m_{h}^{+}u\mathbf m_{h}^{+}v+\frac{h^{2}}{4}\mathbf D_{h}^{+}u\mathbf D_{h}^{+}v,\quad\ \ \ \mathbf D_{h}^{+}(uv)=\mathbf D_{h}^{+}u\mathbf m_{h}^{+}v+\mathbf m_{h}^{+}u\mathbf D_{h}^{+}v,\\ 
\label{2.3}&\Delta_{h}(uv)=\Delta_{h}u{\mathbf m}_{h}v+2{\mathbf D}_{h}u{\mathbf D }_{h}v+{\mathbf m}_{h}u\Delta_{h}v.
\end{align}
\end{lem}

\begin{lem}  Let $ u$ and $v$ be discrete functions defined on $\overline G_h$ such that $ v_{0}=v_{N+1}=0$. Then we have the following identities:
\begin{align} 
&2\int_{G_h}uv\mathbf D_{h}v=-\int_{G_h}\mathbf D_{h}u|v|^{2}+\frac{h^{2}}{2}\int_{G_h^-}\mathbf D_{h}^{+}u|\mathbf D_{h}^{+}v|^{2},\\
\label{2.5}&\int_{G_h}u\Delta_{h}v=-\int_{G_h^{-}}\mathbf D_{h}^{+}u\mathbf D_{h}^{+}v-u_{0}(\mathbf D_{h}^{+}v)_{0}+u_{N+1}(\mathbf D_{h}^{-}v)_{N+1},\\ 
&\int_{G_h}uv\Delta_{h}v=-\int_{G_h^-}\mathbf m_{h}^{+}u|\mathbf D_{h}^{+}v|^{2}+\frac{1}{2}\int_{G_h}\Delta_{h}u|v|^{2},\\ 
&2\int_{G_h}u\mathbf D_{h}v\Delta_{h}v =-\int_{G_h^-}\mathbf D_{h}^{+}u|\mathbf D_{h}^{+}v|^{2}
+u_{N+1}|(\mathbf D_h^{-}v)_{N+1}|^{2}-u_{0}|(\mathbf D_{h}^{+}v)_{0}|^{2}.
\end{align}
\end{lem}

\noindent{\bf Remark 2.1.}\ According to the detailed proof of discrete integration by parts (\ref{2.5}) in [\ref{Bau2013}], we know that it also holds for $v$ without $ v_{0}=v_{N+1}=0$. So we can use (\ref{2.5}) to obtain (\ref{3.60}) in next section.

\subsection{Main results}

 To state our first result, we introduce several weight functions which will be used in our discrete Carleman estimate. Let $\lambda$ and $s$ be two large parameters. For some $x^*<0$, $t_0\in (0,T)$ and a small parameter $\beta>0$,  we define the regular weight functions by 
 \begin{equation}
 	\varphi(x,t)=e^{\lambda \psi(x,t)} ,\quad \theta(x,t)=e^{s\varphi(x,t)},\quad r(x,t)=\theta^{-1}(x,t)
 \end{equation}
 with
 \begin{equation}\label{2.9}
 	\psi(x,t)=|x-x^*|^{2}-\beta |t-t_0|^{2}.
 \end{equation} 

 Moreover, in the following, we will use $C$ to denote generic positive constants depending on $x^*, L, T, \beta$, but independent of $s$ and $\lambda$. Similarly, $C(\lambda)$ denote constants also depending on $\lambda$.  Moreover, we use the notation $\mathcal O_\lambda(\gamma)$, which satisfies $|\mathcal O_\lambda(\gamma)|\leq C(\lambda)|\gamma|$ with a constant $C(\lambda)$. All of these notations may vary from line to line and are independent of $h$. 
 
 According to Proposition 2.9 and Lemma 2.12 in [\ref{Bau2013}], we have the following asymptotic expansion properties: 
 \begin{align}\label{2.10}
 \theta{\mathbf D_{h}r}=-s\lambda A_{1},\quad \theta{\Delta_{h}r}=s^{2}\lambda^{2}A_{2}-s\lambda^{2}A_{3}-s\lambda A_{4}
 \end{align}
 where
\begin{align}\label{2.11}
			A_{j}=f_j+\mathcal{O}_{\lambda}(sh),\quad j=1,2,3,4 
\end{align}
with
\begin{align*}
f_1=\varphi\partial_{x}\psi,\quad f_2=\varphi^{2}|\partial_{x}\psi|^{2},\quad f_{3}=\varphi|\partial_{x}\psi|^{2},\quad f_{4}=\varphi\partial_{xx}\psi.
\end{align*}
Further, we also have for $j=1,2,3,4$ that
\begin{align}\label{2.12}
\left\{
\begin{array}{ll}
\mathbf m_h^{\pm}A_j=f_j+\mathcal{O}_{\lambda}(sh),& \partial_t\mathbf m_h^+A_j=\partial_t f_j+\mathcal O_{\lambda}(sh),\\
\partial_t A_j=\partial_t f_j+\mathcal O_{\lambda}(sh),&\mathbf D_{h} A_j=\partial_{x} f_j+\mathcal O_{\lambda}(sh)=\mathbf D_{h}^{\pm}A_j+\mathcal O_{\lambda}(sh),\\
\Delta_h A_j=\partial_{xx}f_j+\mathcal O_{\lambda}(sh),&\partial_t\mathbf D_{h}^{\pm}A_j=\partial_{xt}f_j+\mathcal O_{\lambda}(sh).
\end{array}
\right.
\end{align}

The first main result in this paper is the following uniform Carleman estimate for the semi-discrete stochastic parabolic equation.

\begin{thm}
 Let $f_{h}\in L^2_{\mathcal F}(0,T; L^2(G_h))$, $g_{h}\in L^2_{\mathcal F}(0,T; H^1(G_h))$. For the parameter $\lambda\geq 1$ sufficiently large, there exist positive constant $C$ depending on $x^*, L, T, \beta$, and positive constants $s_{0}$, $\varepsilon, h_{0}, C(\lambda)$ also depending on $\lambda$, all of which are independent of $h$ such that 
\begin{align}\label{2.13}
		&\mathbb{E} \int_{Q_h}\frac{1}{s\varphi}\theta^{2}|\Delta_{h}y_{h}|^{2}{\rm d}t+ \mathbb{E}\int_{Q_h^-}s\lambda^{2}\varphi\theta^{2}|\partial_{h}^{+}y_{h}|^{2}{\rm d}t+\mathbb{E} \int_{Q_h}s^{3}\lambda^{4}\varphi^{3}\theta^{2}|y_{h}|^{2}{\rm d}t\nonumber\\
&+\mathbb{E}\int_{Q_h^-}s\lambda^{2}\varphi\theta^{2}|g_{h}|^{2}{\rm d}t\nonumber\\
\le&C\mathbb{E}\int_{Q_h}\theta^{2}|f_{h}|^{2}{\rm d}t+C\mathbb{E}\int_{Q_h^-}s\varphi\theta^{2}|\partial_{h}^{+}g_{h}|^{2}{\rm d}t+C\mathbb{E}
\int_{0}^{T}s\lambda\varphi\theta^{2}(x_{N+1},t)|(\mathbf D_{h}^{-}y_{h})_{N+1}|^{2}{\rm d} t\nonumber\\
&+ C(\lambda)s^{2}e^{C(\lambda)s}||y_{h}(T)||_{L^2(\Omega,\mathcal F_T,\mathbb P;L^2(G_h))}^{2}
\end{align}
for all $h\in (0,h_{0} )$, $s\in (s_{0},\sqrt{\varepsilon /h})$ and all $y_{h} \in L_{\mathcal{F} }^{2} (0,T;H^{2}(G_{h}))$ satisfying
\begin{equation}\label{2.14}
	\left\{\begin{array}{ll}
		{\rm d}y_{h}-\Delta _{h}y_{h}{\rm d}t=f_{h}{\rm d}t+g_{h}{\rm d}B(t),&(x_h,t)\in Q_{h}, \\
		y_{h}=0, &(x_h,t)\in\Sigma_{h},\\
        y_{h}(0)=0,&x_h\in G_{h}.
        \end{array}
	\right.
\end{equation}
\end{thm}

\vspace{2mm}

\noindent{\bf Remark 2.2.}\ In comparison with the existing discrete Carleman estimates for parabolic equations [\ref{Boy2020}, \ref{Boy2011}, \ref{Cas2021}], we introduce a regular weight function in our Carleman estimate. The reason is that the solution of a stochastic parabolic equation is not differentiable with respect to time variable. This leads to that the method by applying Carleman estimate to $u_t$ could not directly employed for the inverse random source problem.  
Consequently, we have to choose the regular weight function to put the random source term on the left-side of the Carleman estimate. Unfortunately, there is still a first-order difference term left on the right-hand side of the Carleman estimate. This means that the unknown random source function has to satisfy condition (\ref{2.15}) in the proof of the stability result.

\vspace{2mm}

\noindent{\bf Remark 2.3.}\ In this Carleman estimate, the boundary term at $x_{N+1}$ could be replaced by the one at $x_0$. In fact, if we choose $x^*$ in weight function (\ref{2.9}) such that $x^*>L$, we can obtain $\psi_x<0$ for all $x\in G$. Then the boundary term left in (\ref{3.43}) is  the one at $x_0$.

\vspace{2mm}

\noindent{\bf Remark 2.4.}\ In this Carleman estimate, we use a special function $|x-x^*|^2$  in $\psi$. In fact, we can choose a general form 
\begin{align*}
\psi(x,t)=d(x)-\beta|t-t_0|^2
\end{align*}
with function $d$ such that $|d_x|>0$ in $G$ to guarantee  (\ref{3.39}) and  (\ref{3.40}).

\vspace{2mm}

In order to deal with the discrete Cauchy problem, we need the following  discrete Carleman estimate with non-homogeneous boundary condition.

\begin{thm}
 Let $f_{h}\in L^2_{\mathcal F}(0,T; L^2(G_h))$, $g_{h}\in L^2_{\mathcal F}(0,T; H^1(G_h))$ and $\gamma_1, \gamma_2\in L^2_\mathcal F(\Omega;H^1(0,$ $T))$ satisfying  compatibility condition $\gamma_1(0)=\gamma_2(0)=0$. For the parameter $\lambda\geq 1$ sufficiently large, there exist positive constant $C$ depending on $x^*, L, T, \beta$, and positive constants $s_{0}$, $\varepsilon, h_{0}, C(\lambda)$ also depending on $\lambda$, all of which are independent of $h$ such that
\begin{align}\label{1-2.15}
		&\mathbb{E} \int_{Q_h}\frac{1}{s\varphi}\theta^{2}|\Delta_{h}y_{h}|^{2}{\rm d}t+ \mathbb{E}\int_{Q_h^-}s\lambda^{2}\varphi\theta^{2}|\partial_{h}^{+}y_{h}|^{2}{\rm d}t+\mathbb{E} \int_{Q_h}s^{3}\lambda^{4}\varphi^{3}\theta^{2}|y_{h}|^{2}{\rm d}t\nonumber\\
&+\mathbb{E}\int_{Q_h^-}s\lambda^{2}\varphi\theta^{2}|g_{h}|^{2}{\rm d}t\nonumber\\
\le&C\mathbb{E}\int_{Q_h}\theta^{2}|f_{h}|^{2}{\rm d}t+C\mathbb{E}\int_{Q_h^-}s\varphi\theta^{2}|\partial_{h}^{+}g_{h}|^{2}{\rm d}t+C\mathbb{E}
\int_{0}^{T}s\lambda\varphi\theta^{2}(x_{N+1},t)|(\mathbf D_{h}^{-}y_{h})_{N+1}|^{2}{\rm d} t\nonumber\\
&+C(\lambda)s^3e^{C(\lambda)s}\sum_{i=1}^2\|\gamma_i\|^2_{L^2(\Omega;H^1(0,T))}+ C(\lambda)s^{2}e^{C(\lambda)s}||y_{h}(T)||_{L^2(\Omega,\mathcal F_T,\mathbb P;L^2(G_h))}^{2}
\end{align}
for all $h\in (0,h_{0} )$, $s\in (s_{0},\sqrt{\varepsilon /h})$ and all $y_{h} \in L_{\mathcal{F} }^{2} (0,T;H_h^{2}(G_{h}))$ satisfying
\begin{equation}\label{2.14}
	\left\{\begin{array}{ll}
		{\rm d}y_{h}-\Delta _{h}y_{h}{\rm d}t=f_{h}{\rm d}t+g_{h}{\rm d}B(t),&(x_h,t)\in Q_{h}, \\
        (y_h)_0=\gamma_1,\ (y_h)_{N+1}=\gamma_2,& t\in (0,T),\\
        y_{h}(0)=0,&x_h\in G_{h}.
        \end{array}
	\right.
\end{equation}
\end{thm}

\noindent{\bf Remark 2.5.}\ The duality argument introduced in [\ref{Imanuvilov2001}]  is a main tool to handle non-homogeneous boundary conditions when proving Carleman estimates. It seems that by employing this duality argument, we could obtain a weaker Carleman estimate with $\gamma_1,\gamma_2\in L^2_\mathcal F(\Omega;L^2(0,T))$. In this case, the second derivative term could not be included in the left-hand side of (\ref{1-2.15}). However,  as mentioned in [\ref{Ya2009}], this regularity of $\gamma_1, \gamma_2\in L^2_\mathcal F(\Omega;H^1(0, T))$ is necessary to establish the stability for the Cauchy problem. Therefore, we use a simple method to make the boundary conditions homogeneous, and then use Theorem 2.3 to prove (\ref{1-2.15}).

\vspace{2mm}

Based on the first Carleman estimate, we obtain the uniform stability result with respect to the mesh size $h$ for our discrete inverse random source problem.

\begin{thm} Let  $a_h, b_h\in L^\infty_\mathcal F(0,T;L^\infty(G_h))$ and $g^{(j)}_{h}\in L^2_{\mathcal F}(0,T;H^1(G_h))$ for $j=1,2$ such that
\begin{align}\label{2.15}
 \left|\mathbf D_{h}^{+}\left(g^{(1)}_{h}-g^{(2)}_{h}\right)_{i}\right|\leq C\left|\left(g^{(1)}_{h}-g^{(2)}_{h}\right)_i\right|,\quad i=0,1,2,\cdots,N.
 \end{align} 
 Then there exists a positive constant $C$ depending on   $x^*, L, T$ and $\beta$, but independent of $h$ such that 
\begin{align}\label{2.16}
\left\|g_{h}^{(1)}-g_{h}^{(2)}\right\|_{L^2_{\mathcal F}(0,T;L^2(G_h))}\leq& C\left\|\left(\mathbf D _{h}^{-}y^{(1)}_{h}\right)_{N+1}-\left(\mathbf D _{h}^{-}y^{(2)}_{h}\right)_{N+1}\right\|_{L_{\mathcal F}^{2}(\Omega;L^2(0,T))}\nonumber\\
&+C\left\|y_{h}^{(1)}(T)-y_{h}^{(2)}(T)\right\|_{L^2(\Omega,\mathcal F_T,\mathbb P;L^2(G_h))},
\end{align}
where $y^{(j)}_{h}$ is the solution to (\ref{1.3})  corresponding to $g_{h}^{(j)}$ for $j=1,2$, respectively.
\end{thm}

{\noindent\bf Remark 2.6.}\ A special form of unknown source function is $g_h=r(t)R_h$ with known $R_h$ such that 
\begin{align*}
\left|\left(\mathbf D_h^+R_h\right)_i\right|\leq R_1\quad {\rm and}\quad
|(R_h)_i|\geq R_0>0,\quad {\mathbb P}-a.s.
\end{align*}
with  positive constants $R_0$ and $R_1$.
Then we have
\begin{align*}
\left|\mathbf D_{h}^{+}\left(g^{(1)}_{h}-g^{(2)}_{h}\right)_{i}\right|=\left|r^{(1)}-r^{(2)}\right|\left|\left(\mathbf D_h^+R_h\right)_i\right|\leq \frac{R_1}{R_0}\left|\left(g^{(1)}_{h}-g^{(2)}_{h}\right)_{i}\right|.
\end{align*}
Compared to the continuous inverse problem addressed in [\ref{Wu2020}, \ref{Yuan2021}], a similar condition is imposed on the random source function in the continuous setting to investigate the corresponding inverse problem.

\vspace{2mm}

{\noindent\bf Remark 2.7.}\ There is an additional term depending on the mesh size $h$ in the stability result for discrete inverse problem related to hyperbolic equations [\ref{Bau2015}, \ref{Bau2013}]. This is because of the appearance of the term $h\partial_{h}^{+}\partial_{t} y_{h}$ on the right-hand side of the Carleman estimate, which can not be removed according to [\ref{Bau2013}]. However, for parabolic equations, there is no such an additional term in Carleman estimate (\ref{2.13}). This means that our stability result is uniform with respect to $h$.

\vspace{2mm}

{\noindent\bf Remark 2.8.}\ The uniqueness is a direct result from Theorem 2.5. More precisely, under the same assumptions as in Theorem 2.5 and if
\begin{align*}
	\left(\mathbf D_{h}^{-}y^{(1)}_{h}\right)_{N+1}=\left(\mathbf D_{h}^{-}y^{(2)}_{h}\right)_{N+1}\quad{ and}\quad	y^{(1)}_{h}(T)=y^{(2)}_{h}(T),\quad {\mathbb P}-a.s.,
\end{align*}
then $g_{h}^{(1)}=g_{h}^{(2)}$ in $Q_h$, $\mathbb P-$a.s.

\vspace{2mm}

The last main result is  the stability result for our discrete Cauchy problem, which is an application of the second discrete Carleman estimate.

\begin{thm} Let $a_h, b_h\in L^\infty_\mathcal F(0,T;L^\infty(G_h))$, $c_h\in L^\infty_\mathcal F (0,T; W^{1,\infty}(G_h))$, $\xi\in L^2_\mathcal F(\Omega; H^1$ $(0,T))$ and $\eta\in L^2_\mathcal F(\Omega; L^2(0,T))$. Then for any $\epsilon>0$, there exist positive constants $C$, $h^*$ and $\kappa\in (0,1)$ depending on   $x^*, L, T,\beta$ and $\epsilon$, such that
	\begin{align}\label{2.19}
		&\|y_{h}\|_{L^2_\mathcal F(\epsilon,T-\epsilon;H^2(G_{0,h})}\le CM^\kappa\left(\|\xi\|_{L^2(\Omega;H^1(0,T))}+\|\eta\|_{L^2(\Omega;L^2(0,T))}\right)^{1-\kappa}
	\end{align}
for all $h\in (0,h^*)$ and all $y_{h} \in L_{\mathcal{F} }^{2} (0,T;H^{2}(G_{h}))$ satisfying
\begin{align*}
\|y_{h}\|_{L^2_\mathcal F(0,T;H^2(G_h))}\leq M.
\end{align*}
\end{thm}

{\noindent\bf Remark 2.9.}\ The stability result also holds for discrete Cauchy problem with  the lateral boundary data at $x_0$. Since $\epsilon>0$ is arbitrary,  (\ref{2.19}) immediately implies the
uniqueness of the solution of (\ref{1.4})  in $Q_h$  with $\xi(t)=\eta(t)=0$ for $t\in (0,T)$.

\section{Discrete Carleman estimates}\numberwithin{equation}{section}
 
In this section, we prove two discrete Carleman estimates for semi-discrete stochastic parabolic equation, i.e. Theorem 2.3 and Theorem 2.4. The proof follows as close as possible the ideas presented in the classical continuous setting (see e.g. [\ref{Bar2003}, \ref{Tang2009}]), where  Carleman estimates for stochastic parabolic equations are obtained in the continuous setting.

\subsection{Proof of Theorem 2.3}

 Let $l=s\varphi, \theta=e^{l}, r=\theta^{-1}$ and  $Y_{h}=\theta y_{h}$. By It\^{o} formula, we obtain
\begin{align}\label{3.1}
			\theta {\rm d}y_{h}= {\rm d}Y_{h}-\partial_t l Y_{h}{\rm d}t.
\end{align}
By (\ref{2.1}) and (\ref{2.3}) in Lemma 2.1, we have
\begin{align}\label{3.2}
			\theta\Delta_{h}y_{h}{\rm d}t=&\theta\left( \Delta_{h}r\mathbf m_{h}Y_{h}+2\mathbf D_{h}r\mathbf D_{h}Y_{h}
+\mathbf m_{h}r \Delta_{h}Y_{h}\right){\rm  d}t\nonumber\\	
=&\theta{\Delta_{h}r}\left(Y_{h}+\frac{h^{2}}{4}\Delta_{h}Y_{h}\right){\rm d}t
+2\theta\mathbf D_{h}r\mathbf D_{h}Y_{h}{\rm d}t+\theta\left(r+\frac{h^{2}}{4}\Delta_{h}r\right)\Delta_{h}Y_{h}{\rm d}t\nonumber\\ 
=&{\theta}{\Delta_{h}r}Y_{h}{\rm d}t+\left(1+\frac{h^{2}}{2}{\theta}\Delta_{h}r\right)\Delta_{h}Y_{h}{\rm d}t
+2{\theta}{\mathbf D_{h}r}\mathbf D_{h}Y_{h}{\rm d}t.
\end{align}
Then from (\ref{3.1}) and (\ref{3.2}), it follows
\begin{align}\label{3.3}
			&\theta \left({\rm d}y_{h}-\Delta_{h}y_{h}{\rm d}t\right)\nonumber\\
=&{\rm d}Y_{h}-\partial_tl Y_{h}{\rm d}t-{\theta}{\Delta_{h}r}Y_{h}{\rm d}t-\left(1+\frac{h^{2}}{2}{\theta}\Delta_{h}r\right)\Delta_{h}Y_{h}{\rm d}t
-2{\theta}{\mathbf D_{h}r}\mathbf D_{h}Y_{h}{\rm d}t.
\end{align}
According to (\ref{2.10}), we rewrite (\ref{3.3}) as the following form:
\begin{equation}\label{3.4}
	\theta ({\rm d}y_{h}-\Delta_{h}y_{h}{\rm d}t)=I_{1} +I{\rm d}t,
\end{equation}
where
$$
\begin{aligned}
		        &I_{1}={\rm d}Y_{h}+2s\lambda A_{1}\mathbf D_{h}Y_{h}{\rm d}t+\Psi Y_{h}{\rm d}t,\\
&I=-(1+A_0)\Delta_{h}Y_{h}-s^{2}\lambda ^{2}A_{2}Y_{h}+(-\partial_t l +s\lambda ^{2} A_{3} +s\lambda A_{4}-\Psi )Y_{h}.
\end{aligned}
$$
with
\begin{align*}
A_0=\frac{h^2}{2}\left(s^2\lambda^2A_2-s\lambda^2A_3-s\lambda A_4\right),\quad\Psi=\tau\partial_{xx} l.
\end{align*}
Here $\tau$ is a  positive constant such that $\tau\in(3/2,3)$. 
By (\ref{2.12}) we have the following properties of $A_0$: 
\begin{align}
A_0=\frac{h^2}{2}\left(s^2\lambda^2f_2-s\lambda^2f_3-s\lambda f_4+s^2\mathcal O_\lambda(sh)\right)=\mathcal O_{\lambda}(sh).
\end{align}
Similarly,
\begin{align}\label{3.6}
\mathbf m_h^{+}A_0=\mathcal O_{\lambda}(sh),\quad \mathbf D_h^{+} A_0=\mathcal O_{\lambda}(sh),\quad \Delta_h A_0=\mathcal O_{\lambda}(sh).
\end{align}
Furthermore, multiplying $I$ in both sides of (\ref{3.4}), integrating the result equality over $Q_h$ and taking mathematical expectation, we find
\begin{align}\label{3.7}
		\mathbb E\int_{Q_h}\theta I({\rm d} y_{h}-\Delta_{h}y_{h}{\rm d}t)=\mathbb E\int_{Q_h}II_{1} +\mathbb E\int_{Q_h}I^2{\rm d}t.
\end{align}

Next, we provide a detailed calculation of the term involving $II_1$ and then find positive lower bounds for the terms related to   ${\rm d}Y_h$ and ${\rm d}t$, respectively. Finally, we combine these estimates to complete our proof. 
To do this, we divide the subsequent proof into the following several steps. For clarity, we first split the term of $II_{1}$ into a sum of nine terms
\begin{equation}
	\mathbb E\int_{Q_h}II_{1}=\sum_{i,j=1}^{3} I_{ij},
\end{equation}
 where $I_{ij}$ is the inner product of $i$-th term of $I_1$ and the $j$-th term of $I$.

 \subsubsection{Estimates involving the terms of ${\rm d}Y_h$}

 By  discrete identities (\ref{2.1}) and discrete integration by parts (\ref{2.5}), together with $$\left({\rm d}Y_h\right)_0=\left({\rm d}Y_h\right)_{N+1}=0,\quad t\in [0,T]$$ 
 due to $y_h=0$ on $\Sigma_h$, we obtain
\begin{align}
	I_{11}=&-\mathbb{E} \int_{Q_h}(1+A_0)\Delta _{h}Y_{h}{\rm d}Y_{h}\nonumber\\ 
      =&\mathbb{E}\int_{Q_h^{-}}\mathbf D_{h}^{+}Y_h\mathbf D_{h}^{+}({\rm d}Y_{h}+A_0{\rm d}Y_h)+\mathbb{E}\int_{0}^{T}(\mathbf D_{h}^{+}Y_{h})_{0}({\rm d}Y_{h}+A_0{\rm d}Y_h)_{0}\nonumber\\
&-\mathbb{E}\int_{0}^{T}(\mathbf D_{h}^{-}Y_{h})_{N+1}({\rm d}Y_{h}+A_0{\rm d}Y_h)_{N+1}\nonumber\\
	=&\mathbb E\int_{Q_h^-}\left(1+\frac{h}{2}\mathbf D_h^{+}A_0+\mathbf m_h^{+}A_0\right)\mathbf D_h^{+}Y_h\mathbf D_h^{+}{\rm d}Y_h+\mathbb E\int_{Q_h^-}\mathbf D_h^+A_0\mathbf D_h^+Y_h{\rm d}Y_h.
\end{align}
By using It\^{o} formula, we further obtain
\begin{align}
I_{11}=&X_1+Y_1+Z_1+\mathbb E\int_{Q_h^-}B_1|\mathbf D_h^+ Y_h|^2{\rm d}t,
\end{align}
where
\begin{align*}
&X_{1}=\frac{1}{2}\mathbb E\int_{Q_h^-}{\rm d}\left(\left(1+\frac{h}{2}\mathbf D_h^{+}A_0+\mathbf m_h^{+}A_0\right)|\mathbf D_h^+Y_h|^2\right),\\
&Y_1=-\frac{1}{2}\mathbb E\int_{Q_h^-}\left(1+\frac{h}{2}\mathbf D_h^{+}A_0+\mathbf m_h^{+}A_0\right)|\mathbf D_h^+{\rm d} Y_h|^2,\\
&Z_1=\mathbb E\int_{Q_h^-}\mathbf D_h^+A_0\mathbf D_h^+Y_h{\rm d}Y_h,\\
&B_1=-\frac{1}{2}\partial_t\left(1+\frac{h}{2}\mathbf D_h^{+}A_0+\mathbf m_h^{+}A_0\right).
\end{align*}
Using It\^{o} formula again, we obtain
\begin{align}
		I_{12}=&-\mathbb{E}\int_{Q_h}s^{2}\lambda ^{2}  A_{2} Y_{h}{\rm d}Y_{h}=X_2+Y_2+\mathbb{E}\int_{Q_h}D_1|Y_h|^2{\rm d}t,
\end{align}
where
\begin{align*}
X_2=&-\frac{1}{2}\mathbb{E}\int_{Q_h}{\rm d}(s^{2}\lambda ^{2}A_2|Y_{h}| ^{2}),\\
Y_2=&\frac{1}{2}\mathbb{E}\int_{Q_h}s^{2}\lambda ^{2}A_{2}|{\rm d}Y_{h}| ^{2},\\
D_1=&\frac{1}{2}s^{2}\lambda ^{2}\partial_tA_{2}.
\end{align*}
Similarly,
\begin{align}
	I_{13}=&\mathbb{E} \int_{Q_h}(-\partial_t l+s\lambda ^{2}A_{3}+s\lambda A_{4}-\Psi )Y_{h}{\rm d}Y_{h}=X_3+Y_3+\mathbb{E} \int_{Q_h}D_2|Y_{h}|^{2}{\rm d}t,
\end{align}
where
\begin{align*}
X_3=&\frac{1}{2} \mathbb{E} \int_{Q_h}{\rm d}\left(\left(-\partial_t l +s\lambda ^{2}A_{3}+s\lambda A_{4}-\Psi \right)|Y_{h}|^{2}\right),\\
	Y_3=&-\frac{1}{2} \mathbb{E} \int_{Q_h}(-\partial_t l +s\lambda ^{2}A_{3}+s\lambda A_{4}-\Psi )|{\rm d}Y_{h}|^{2},\\
D_2=&-\frac{1}{2} (-\partial_{tt} l+s\lambda ^{2}\partial_tA_{3}+s\lambda \partial_t A_{4}-\partial_{t}\Psi ).
\end{align*}
Therefore, we have
\begin{align}\label{3.13}
\sum_{j=1}^3I_{1j}=\sum_{i=1}^{3} X_i+\sum_{i=1}^3Y_i+Z_1+\sum_{i=1}^2\mathbb{E} \int_{Q_h}D_i|Y_{h}|^{2}{\rm d}t+\mathbb E\int_{Q_h^-}B_1|\mathbf D_h^+ Y_h|^2{\rm d}t.
\end{align}

Using (\ref{3.6}),  we obtain
\begin{align}
&1+\frac{h}{2}\mathbf D_h^{+}A_0+\mathbf m_h^{+}A_0=1+\mathcal O_{\lambda}(sh).
\end{align}
If we choose $\varepsilon=\varepsilon(\lambda)$ sufficiently small such that $|\mathcal O_\lambda(sh)|\leq \frac{1}{2}$, we further have
\begin{align}\label{3.15}
1+\frac{h}{2}\mathbf D_h^{+}A_0+\mathbf m_h^{+}A_0\geq \frac{1}{2}.
\end{align}
On the other hand, we have
\begin{align}\label{3.16}
&-\partial_tl+s\lambda ^{2}A_{3}+s\lambda A_{4}-\Psi=-s\partial_t\varphi+s\lambda ^{2}A_{3}+s\lambda A_{4}-\tau\partial_{xx} l\nonumber\\
=&-s\lambda\varphi\psi_t+(1-\tau)(s\lambda^2\varphi|\partial_x\psi|^2
+s\lambda\varphi\partial_{xx}\psi)+s\mathcal O_\lambda(sh).
\end{align}
 By (\ref{3.15}) and (\ref{3.16}), together with $Y_{h}(0)=0$ due to $y_h(0)=0$, we then use integration by parts with respect to time $t$ to yield
\begin{align*}
		\sum_{i=1}^{3} X_i=&\left.\frac{1}{2}\mathbb E\int_{G_h^-}\left(\left(1+\frac{h}{2}\mathbf D_h^{+}A_0+\mathbf m_h^{+}A_0\right)|\mathbf D_h^+Y_h|^2\right)\right|_{t=T}\nonumber\\
&-\left.\frac{1}{2}\mathbb{E}\int_{G_h}(s^{2}\lambda ^{2}A_2|Y_{h}| ^{2})\right|_{t=T}\nonumber\\
&+\left.\frac{1}{2} \mathbb{E} \int_{G_h}\left(\left(-\partial_tl +s\lambda ^{2}A_{3}+s\lambda A_{4}-\Psi \right)|Y_{h}|^{2}\right)\right|_{t=T}
\end{align*}
\begin{align}\label{3.17}
&\hspace{1.5cm}\geq\left.\frac{1}{4}\mathbb E\int_{G_h^-}\left(|\mathbf D_h^+Y_h|^2\right)\right|_{t=T}-C(\lambda)s^2(1+|\mathcal O_{\lambda}(sh)|)\left.\mathbb{E}\int_{G_h}(|Y_{h}| ^{2})\right|_{t=T}\nonumber\\
&\hspace{1.8cm}-C(\lambda)s\left(1+|\mathcal O_\lambda(sh)|\right)\left.\mathbb{E} \int_{G_h}\left(|Y_{h}|^{2}\right)\right|_{t=T}\nonumber\\
&\hspace{1.5cm}\geq -C(\lambda ) s^2\left.\mathbb{E}\int_{G_h} \left(\theta^2|y_{h}|^{2}\right)\right|_{t=T}.
\end{align}

Next we estimate the terms involving $Y_i$ in (\ref{3.13}). By (\ref{3.6}), we obtain
\begin{align}\label{3.18}
\sum_{i=1}^{3} Y_i=&-\frac{1}{2}\mathbb E\int_{Q_h^-}(1+\mathcal O_{\lambda}(sh))|\mathbf D_h^+{\rm d} Y_h|^2\nonumber\\
&+\frac{1}{2}\mathbb{E}\int_{Q_h}\left(s^{2}\lambda ^{2}A_2+\partial_t l-s\lambda ^{2}A_{3}-s\lambda A_{4}+\Psi\right)|{\rm d}Y_{h}| ^{2}.
\end{align}
Obviously,
\begin{align}
		|\mathbf D_{h}^{+} {\rm d} Y_{h}|^{2}
		=&|\mathbf D_{h}^{+}(\theta {\rm d}y_{h})|^{2}=|\mathbf D_{h}^{+}\theta \mathbf m_{h}^{+}{\rm d}y_{h}+\mathbf m_{h}^{+}\theta \mathbf D_{h}^{+}{\rm d}y_{h}|^{2}\nonumber\\
		=&|\mathbf D_{h}^{+}\theta {\rm d} y_{h}+h\mathbf D_{h}^{+}\theta \mathbf D_{h}^{+}{\rm d}y_h+\theta \mathbf D_{h}^{+}{\rm d}y_{h}|^{2}\nonumber\\
		=&|\mathbf D_{h}^{+}\theta |^{2}|{\rm d}y_{h}|^{2}+\left(h\mathbf D_{h}^{+}\theta+\theta\right)^2|\mathbf D_{h}^{+}{\rm d}y_{h}|^{2}+2\left(h|\mathbf D_{h}^{+}\theta |^{2}+\theta\mathbf D_{h}^{+}\theta\right){\rm d}y_{h}\mathbf D_{h}^{+}{\rm d}y_{h}.
\end{align}
By using Young's inequality, we further obtain
\begin{align}\label{3.20}
		|\mathbf D_{h}^{+} {\rm d}Y_{h}|^{2}
		\le&|\mathbf D_{h}^{+}\theta |^{2}|{\rm d}y_{h}|^{2}+2h^{2}|\mathbf D_{h}^{+}\theta |^{2}|\mathbf D_{h}^{+}{\rm d}y_{h}|^{2}+2\theta ^{2}|\mathbf D_{h}^{+}{\rm d}y_{h}|^{2}\nonumber\\
    &+\frac{1}{2}h^{2}\partial_{xx} l|\mathbf D_{h}^{+}\theta|^{2}|{\rm d} y_{h}|^{2}+\frac{2}{\partial_{xx} l}|{\mathbf D}_{h}^{+}\theta|^{2}|\mathbf D_{h}^{+}{\rm d}y_{h}|^{2}\nonumber\\
		&+\frac{1}{2}\partial_{xx} l\theta^{2}|{\rm d}y_{h}|^{2}+\frac{2}{\partial_{xx} l}|\mathbf D_{h}^{+}\theta|^{2} |\mathbf D_{h}^{+}{\rm d}y_{h}|^{2}\nonumber\\ 
\le&\left(|\mathbf D_{h}^{+}\theta |^{2}+\frac{1}{2} h^{2}\partial_{xx} l|\mathbf D_{h}^{+}\theta|^{2}+\frac{1}{2} \partial_{xx} l\theta^{2}\right)|{\rm d}y_{h}|^{2}\nonumber\\
		&+\left(2h^{2}|\mathbf D_{h}^{+}\theta|^{2}+2\theta ^{2}+\frac{4}{\partial_{xx} l}|\mathbf D_{h}^{+}\theta|^{2}\right)|\mathbf D_{h}^{+}{\rm d}y_{h}|^{2}.
\end{align}
By Taylor formula, we have
\begin{align}\label{3.21}\left(\mathbf D_{h}^{+}\theta\right)_j=\int_0^1\partial_x\theta(x_j+\sigma h,t){\rm d}\sigma=\left(\partial_{x}\theta\right)_{j}+s\mathcal{O}_{\lambda}(sh)\theta_{j}.
\end{align} Similarly,
\begin{align}\label{3.22} &\left(\mathbf \Delta_h\theta\right)_j=\left(\partial_{xx}\theta\right)_j +s^2\mathcal{O}_{\lambda}(sh)\theta_j.
\end{align}
Therefore, we obtain
\begin{align}\label{3.23}
&|\mathbf D_{h}^{+}\theta |^{2}+\frac{1}{2} h^{2}\partial_{xx} l|\mathbf D_{h}^{+}\theta|^{2}+\frac{1}{2} \partial_{xx} l\theta^{2}\nonumber\\
\leq &\left(1+\frac{1}{2}h^2(s\lambda^{2}\varphi|\partial_{x}\psi|^{2}+s\lambda\varphi\partial_{xx}\psi)\right)
\left(s^{2}\lambda^{2}\varphi^{2}|\partial_{x}\psi|^{2}+s^2\left|\mathcal{O}_{\lambda }(sh)\right|\right)\theta^{2}\nonumber\\
&+\frac{1}{2}\left(s\lambda^{2}\varphi|\partial_{x}\psi|^{2}+s\lambda\varphi\partial_{xx}\psi\right)\theta^2\nonumber\\
\leq & \left(s^{2}\lambda^{2}\varphi^{2}|\partial_{x}\psi|^{2}+\frac{1}{2}s\lambda^{2}\varphi|\partial_{x}\psi|^{2}
+\frac{1}{2}s\lambda\varphi\partial_{xx}\psi+s^2\left|\mathcal{O}_{\lambda }(sh)\right|\right)\theta^{2}
\end{align}
and
\begin{align}\label{3.24}
2h^{2}|\mathbf D_{h}^{+}\theta|^{2}+2\theta ^{2}+\frac{4}{\partial_{xx} l}|\mathbf D_{h}^{+}\theta|^{2}\leq & 2h^2\left(s^{2}\lambda^{2}\varphi^{2}|\partial_{x}\psi|^{2}+s^2\left|\mathcal{O}_{\lambda}(sh)\right|\right)\theta^{2}
+2\theta^2\nonumber\\
&+\frac{4\left(s^{2}\lambda^{2}\varphi^{2}|\partial_{x}\psi|^{2}+s^2\left|\mathcal{O}_{\lambda}(sh)\right|\right)\theta^2}
{s\lambda^{2}\varphi|\partial_{x}\psi|^{2}+s\lambda\varphi\partial_{xx}\psi}\nonumber\\
\leq & \left(2+4s\varphi+s|\mathcal O_\lambda(sh)|\right)\theta^2.
\end{align}
Combining (\ref{3.23}), (\ref{3.24}) and (\ref{3.20}) we obtain
\begin{align}\label{3.25}
		|\mathbf D_{h}^{+} {\rm d}Y_{h}|^{2}\le&\left(s^{2}\lambda^{2}\varphi^{2}|\partial_{x}\psi|^{2}+\frac{1}{2}s\lambda^{2}\varphi|\partial_{x}\psi|^{2}
+\frac{1}{2}s\lambda\varphi\partial_{xx}\psi+s^2\left|\mathcal{O}_{\lambda }(sh)\right|\right)\theta^{2}|{\rm d}y_h|^2\nonumber\\
&+C\left(s\varphi+s|\mathcal{O}_{\lambda }(sh)|\right)\theta^{2}|\mathbf D_{h}^{+}{\rm d}y_h|^{2}.
\end{align}
On the other hand, from (\ref{2.11}) and (\ref{3.16}), it follows that
\begin{align}\label{3.26}
&s^{2}\lambda ^{2}A_2+\partial_t l-s\lambda ^{2}A_{3}-s\lambda A_{4}+\Psi\nonumber\\
=&s^2\lambda^2\varphi^2|\partial_x\psi|^2+s\lambda\varphi\psi_t+(\tau-1)(s\lambda^2\varphi|\partial_x\psi|^2
+s\lambda\varphi\partial_{xx}\psi)+s\mathcal O_\lambda(sh).
\end{align}
Then, noticing that ${\rm d}y_h=0$ at $x=x_0$ due to $\left(y_h\right)_0=0$ and substituting (\ref{3.25}), (\ref{3.26}) into (\ref{3.18}), we obtain the following estimate
\begin{align}
\sum_{i=1}^{3} Y_i\geq &-\frac{1}{2}\mathbb E\int_{Q_h^-}(s^{2}\lambda^{2}\varphi^{2}|\partial_{x}\psi|^{2}+\frac{1}{2}s\lambda^{2}\varphi|\partial_{x}\psi|^{2}
+\frac{1}{2}s\lambda\varphi\partial_{xx}\psi)\theta^{2}|{\rm d}y_h|^2\nonumber\\
&-s^2\mathbb E\int_{Q_h^-} |\mathcal O_\lambda(sh)|\theta^{2}|{\rm d}y_h|^2-C\mathbb E\int_{Q_h^-}\left(s\varphi+s|\mathcal{O}_{\lambda }(sh)|\right)\theta^{2}|\mathbf D_{h}^{+}{\rm d}y_h|^{2}\nonumber\\
&+\frac{1}{2}\mathbb{E}\int_{Q_h^-}\left(s^2\lambda^2\varphi^2|\partial_x\psi|^2+s\lambda\varphi\psi_t
+(\tau-1)(s\lambda^2\varphi|\partial_x\psi|^2
+s\lambda\varphi\partial_{xx}\psi)\right)\theta^{2}|{\rm d}y_h|^2\nonumber\\
\geq&\frac{1}{4}\mathbb{E}\int_{Q_h^-}\left(\tau-\frac{3}{2}\right)s\lambda^2\varphi|\partial_x\psi|^2\theta^{2}|{\rm d}y_h|^2-s^2\mathbb E\int_{Q_h^-} |\mathcal O_\lambda(sh)|\theta^{2}|{\rm d}y_h|^2\nonumber\\
&-C\mathbb E\int_{Q_h^-}\left(s\varphi+s|\mathcal{O}_{\lambda }(sh)|\right)\theta^{2}|\mathbf D_{h}^{+}{\rm d}y_h|^{2}.
\end{align}
where we choose $\lambda$ sufficiently small such that $$\frac{1}{4}\left(\tau-\frac{3}{2}\right)\lambda|\psi_x|^2\geq \frac{1}{2}\left(\tau-\frac{3}{2}\right)|\partial_{xx}\psi|+\frac{1}{2}|\psi_t|,\quad t\in [0,T].$$
Since $s\in (s_{0},\sqrt{\varepsilon /h})$, there exists $\varepsilon(\lambda)$ sufficiently small such that
\begin{align}
s^2|\mathcal O_{\lambda}(sh)|\leq \frac{\varepsilon}{h}C(\lambda)sh\leq \frac{1}{8}s.
\end{align}
Thus, together with $|{\rm d}y_h|^2=|g_h|^2{\rm d}t$, we obtain
\begin{align}\label{3.29}
\sum_{i=1}^{3} Y_i\geq&\frac{1}{4}\mathbb{E}\int_{Q_h^-}\left(\tau-2\right)s\lambda^2\varphi|\partial_x\psi|^2\theta^{2}|g_h|^2{\rm d}t-C\mathbb E\int_{Q_h^-}s\varphi\theta^{2}|\mathbf D_{h}^{+}g_h|^{2}{\rm d}t.
\end{align}

Now we deal with the cross term $Z_1$ in (\ref{3.13}). By (\ref{2.1}) and It\^{o}'s inequality, we obtain
\begin{align}\label{3.30}
		 &\mathbf D_h^+A_0\mathbf D_{h}^{+}Y_{h}{\rm d}Y_{h}\nonumber\\
=&\mathbf D_h^+A_0(\mathbf D_{h}^{+}\theta \mathbf m_{h}^{+}y_{h}+\mathbf m_{h}^{+}\theta \mathbf D_{h}^{+}y_{h})(\theta {\rm d}y_{h}+\partial_t\theta y_{h}{\rm d}t)\nonumber\\
=&\mathbf D_h^+A_0\left(\mathbf D_{h}^{+}\theta y_{h}+\theta\mathbf D_{h}^{+}y_h+h\mathbf D_{h}^{+}\theta  \mathbf D_{h}^{+}y_{h}\right)(\theta {\rm d}y_{h}+\partial_{t}\theta y_{h}{\rm d}t)\nonumber\\
=&\frac{1}{2}{\rm d}\left(\mathbf D_h^+A_0\theta \mathbf D_h^+\theta|y_h|^2\right)-\frac{1}{2}\mathbf D_h^+A_0\theta \mathbf D_h^+\theta|{\rm d}y_h|^2-\frac{1}{2}\partial_t(\mathbf D_h^+A_0\theta\mathbf D_h^+\theta)|y_h|^2{\rm d}t\nonumber\\
&+\mathbf D_h^+A_0\partial_t\theta\mathbf D_h^+\theta|y_h|^2{\rm d}t+\mathbf D_h^+A_0\left(\theta+h\mathbf D_h^+\theta\right)\partial_t\theta y_h\mathbf D_h^+y_h{\rm d}t\nonumber\\
&+\mathbf D_h^+A_0\left(\theta+h\mathbf D_h^+\theta\right)\theta \mathbf D_h^+y_h{\rm d}y_h.
\end{align}
By (\ref{3.6}) and (\ref{3.21}), we have
\begin{align}\label{3.31}
\left\{\begin{array}{l}
\mathbf D_h^+A_0\theta \mathbf D_h^+\theta=\mathcal O_{\lambda}(sh)(s\lambda\varphi\partial_x\psi+s\mathcal O_{\lambda}(sh))\theta^2=s\mathcal O_{\lambda}(sh)\theta^2,\\
\mathbf D_h^+A_0\left(\theta+h\mathbf D_h^+\theta\right)\theta=\mathcal O_{\lambda}(sh)\left(1+h(s\lambda\varphi\partial_x\psi+s\mathcal O_{\lambda}(sh))\right)\theta^2=\mathcal O_{\lambda}(sh)\theta^2.
\end{array}
\right.
\end{align}
Then, using (\ref{3.30}) and (\ref{3.31}) we obtain the following estimate for $Z_1$:
\begin{align}\label{3.32}
Z_1\geq&-\left.\mathbb E\int_{G_h^-}s|\mathcal O_{\lambda}(sh)|\left(\theta^2|y_h|^2\right)\right|_{t=T}-\mathbb E\int_{Q_h^-}s|\mathcal O_{\lambda}(sh)|\theta^2|g_h|^2{\rm d}t+\mathbb E\int_{Q_h^-}D_3|Y_h|^2{\rm d}t\nonumber\\
&+\mathbb E\int_{Q_h^-}Ky_h\mathbf D_h^+y_h{\rm d}t+\mathbb E\int_{Q_h^-}\mathcal O_{\lambda}(sh)\theta^2 \mathbf D_h^+y_h{\rm d}y_h,
\end{align}
where
\begin{align*}
&D_3=\left(-\frac{1}{2}\partial_t(\mathbf D_h^+A_0\theta\mathbf D_h^+\theta)+\mathbf D_h^+A_0\partial_t\theta\mathbf D_h^+\theta\right)\theta^{-2},\nonumber\\
&K=\mathbf D_h^+A_0\left(\theta+h\mathbf D_h^+\theta\right)\partial_t\theta .
\end{align*}
By the equation of $y_h$ in (\ref{2.14}), we have 
\begin{align}\label{3.33}\mathbf D_{h}^{+}y_{h}{\rm d}y_{h}=\mathbf D_{h}^{+}y_{h}\Delta _{h}y_{h}{\rm d}t+f_{h}\mathbf D_{h}^{+}y_{h}{\rm d}t+g_{h}\mathbf D_{h}^{+}y_{h}{\rm d}B(t).
\end{align} 
Then,  using (\ref{3.33}) and $$\mathbb E\int_{Q_h^-}\mathcal O_{\lambda}(sh)\theta^2 \mathbf D_h^+y_h{\rm d}y_h=\mathbb E\int_{Q_h}\mathcal O_{\lambda}(sh)\theta^2 \mathbf D_h^+y_h{\rm d}y_h$$ duet to $y_h=0$ at $x=x_0$,  we obtain the following estimate for the last term on the right-hand side of (\ref{3.32}):
\begin{align*}
		&\mathbb E\int_{Q_h^-}\mathcal O_{\lambda}(sh)\theta^2 \mathbf D_h^+y_h{\rm d}y_h\nonumber\\
=&\mathbb E\int_{Q_h}\mathcal O_{\lambda}(sh)\theta^2 \left(\mathbf D_{h}^{+}y_{h}\Delta _{h}y_{h}{\rm d}t+f_{h}\mathbf D_{h}^{+}y_{h}{\rm d}t+g_{h}\mathbf D_{h}^{+}y_{h}{\rm d}B(t)\right)\nonumber\\
=&\mathbb E\int_{Q_h}\mathcal O_{\lambda}(sh)\theta^2 \mathbf D_{h}^{+}y_{h}\Delta _{h}y_{h}{\rm d}t+\mathbb E\int_{Q_h}\mathcal O_{\lambda}(sh)\theta^2 f_{h}\mathbf D_{h}^{+}y_{h}{\rm d}t,
\end{align*}
where we have used
\begin{align*} \mathbb E\int_{Q_h}\mathcal O_{\lambda}(sh)\theta^2 g_{h}\mathbf D_{h}^{+}y_{h}{\rm d}B(t)=0.
\end{align*}
Young's inequality with $\epsilon$ further yields that
\begin{align}\label{3.34}
		&\mathbb E\int_{Q_h^-}\mathcal O_{\lambda}(sh)\theta^2 \mathbf D_h^+y_h{\rm d}y_h\nonumber\\
\geq &-\epsilon\mathbb E\int_{Q_h}\frac{1}{s\varphi}\theta^2 |\Delta _{h}y_{h}|^2{\rm d}t-C(\epsilon)\mathbb E\int_{Q_h}s|\mathcal O_{\lambda}(sh)|\varphi\theta^2 |\mathbf D_{h}^{+}y_{h}|^2{\rm d}t\nonumber\\
&-\mathbb E\int_{Q_h}\theta^2 |f_{h}|^2{\rm d}t,
\end{align}

Therefore, we deduce from  (\ref{3.34}) and (\ref{3.32}) that
\begin{align}\label{3.35}
Z_1\geq&-\left.\mathbb E\int_{G_h^-}s|\mathcal O_{\lambda}(sh)|\left(\theta^2|y_h|^2\right)\right|_{t=T}-\mathbb E\int_{Q_h^-}s|\mathcal O_{\lambda}(sh)|\theta^2|g_h|^2{\rm d}t+\mathbb E\int_{Q_h^-}D_3|y_h|^2{\rm d}t\nonumber\\
&-\epsilon\mathbb E\int_{Q_h}\frac{1}{s\varphi}\theta^2 |\Delta _{h}y_{h}|^2{\rm d}t-C(\epsilon)\mathbb E\int_{Q_h}s|\mathcal O_{\lambda}(sh)|\varphi\theta^2 |\mathbf D_{h}^{+}y_{h}|^2{\rm d}t-\mathbb E\int_{Q_h}\theta^2 |f_{h}|^2{\rm d}t\nonumber\\
&+\mathbb E\int_{Q_h^-}Ky_h\mathbf D_h^+y_h{\rm d}t.
\end{align}

Substituting (\ref{3.17}), (\ref{3.29}) and (\ref{3.35}) into (\ref{3.13}), and noticing that $s|\mathcal O_\lambda(sh)|\leq 1/8$ we obtain
\begin{align}\label{3.36}
\sum_{j=1}^3I_{1j}\geq&\frac{1}{4}\mathbb{E}\int_{Q_h^-}\left(\tau-2\right)s\lambda^2\varphi|\partial_x\psi|^2\theta^{2}|g_h|^2{\rm d}t-C\mathbb E\int_{Q_h^-}s\varphi\theta^{2}|\mathbf D_{h}^{+}g_h|^{2}{\rm d}t\nonumber\\
&-\mathbb E\int_{Q_h}\theta^2 |f_{h}|^2{\rm d}t-C(\lambda)s^{2}e^{C(\lambda)s}||y_{h}(T)||_{L^2(\Omega,\mathcal F_T,\mathbb P;L^2(G_h))}^{2}\nonumber\\
&-\epsilon\mathbb E\int_{Q_h}\frac{1}{s\varphi}\theta^2 |\Delta _{h}y_{h}|^2{\rm d}t-C(\epsilon)\mathbb E\int_{Q_h}\varphi\theta^2 |\mathbf D_{h}^{+}y_{h}|^2{\rm d}t+\sum_{i=1}^3\mathbb{E} \int_{Q_h}D_i|Y_{h}|^{2}{\rm d}t\nonumber\\
&+\mathbb E\int_{Q_h^-}B_1|\mathbf D_h^+ Y_h|^2{\rm d}t+\mathbb E\int_{Q_h^-}Ky_h\mathbf D_h^+y_h{\rm d}t.
\end{align}

 \subsubsection{Estimates involving the terms of ${\rm d}t$}
 
 Proceeding as done in [\ref{Bau2013}] or [\ref{Boy2014}], we apply discrete integrations by parts in Lemma 2.2 to yield
 \begin{flalign*}
	\bullet\quad
	I_{21}=&-2\mathbb{E} \int_{Q_h}s\lambda A_{1}(1+A_0)\mathbf D_{h}Y_{h}\Delta_{h}Y_{h}{\rm d}t\\
=& \mathbb{E} \int_{Q_h^-}s\lambda\mathbf D _{h}^{+}(A_{1}(1+A_0))|\mathbf D _{h}^{+}Y_{h}|^{2}{\rm d}t+ \mathbb{E} \int_{0}^{T}s\lambda (A_{1}(1+A_0))(x_0,t)|(\mathbf D _{h}^{+}Y_{h})_{0}|^{2}{\rm d}t\nonumber\\
&-\mathbb{E} \int_{0}^{T}s\lambda (A_{1}(1+A_0)(x_{N+1},t))|(\mathbf D_{h}^{-}Y_{h})_{N+1}|^{2}{\rm d}t,&
\end{flalign*}
\begin{flalign*}
	\bullet\quad
	I_{22}=&-2\mathbb{E} \int_{Q_h}s^{3}\lambda ^{3}A_{1}A_{2}Y_{h}\mathbf D_{h}Y_{h}{\rm d}t
	\\=&\mathbb{E} \int_{Q_h}s^{3}\lambda ^{3}\mathbf D_{h}(A_{1}A_{2})|Y_{h}|^{2}{\rm d}t-\frac{h^{2}}{2}\mathbb{E} \int_{Q_h^-}s^{3}\lambda ^{3}\mathbf D_{h}^{+}(A_{1}A_{2})|\mathbf D_{h}^{+}Y_{h}|^{2}{\rm d}t,&
\end{flalign*}
\begin{flalign*}
	\bullet\quad
	I_{23}=&2\mathbb{E} \int_{Q_h}s\lambda A_{1}(-\partial_tl+s\lambda ^{2}A_{3}+s\lambda A_{4}-\Psi)Y_{h}\mathbf D _{h}Y_{h}{\rm d}t\\
=&-\mathbb E\int_{Q_h}s\lambda\mathbf D_h\left(A_{1}(-\partial_tl+s\lambda ^{2}A_{3}+s\lambda A_{4}-\Psi)\right)|Y_h|^2{\rm d}t\\
&+\frac{h^2}{2}\mathbb E\int_{Q_h^-}s\lambda\mathbf D_h^+\left(A_{1}(-\partial_tl+s\lambda ^{2}A_{3}+s\lambda A_{4}-\Psi)\right)|\mathbf D_h^+ Y_h|^2{\rm d}t,&
\end{flalign*}
\begin{flalign*}
	\bullet\quad
	I_{31}=&-\mathbb{E} \int_{Q_h} (1+A_0)\Psi Y_{h}\Delta _{h}Y_{h}{\rm d}t\\
  =&\mathbb{E} \int_{Q_h^-}\mathbf m_{h}^{+}\left((1+A_0)\Psi\right)|\mathbf D _{h}^{+}Y_{h}|^{2}{\rm d}t -\frac{1}{2} \mathbb{E} \int_{Q_h}\Delta _{h}((1+A_0)\Psi)|Y_{h}|^{2}{\rm d}t,&
\end{flalign*}
\begin{flalign*}
	\bullet\quad
	I_{32}=&-\mathbb{E} \int_{Q_h}s ^{2}\lambda ^{2}A_{2}\Psi |Y_{h}|^{2}{\rm d}t,&
\end{flalign*}
\begin{flalign*}
	\bullet\quad
	I_{33}=&\mathbb{E} \int_{Q_h}(-\partial_tl+s\lambda ^{2}A_{3}+s\lambda A_{4}-\Psi)\Psi|Y_{h}|^{2}{\rm d}t.&
\end{flalign*}
Then we find that
\begin{align}\label{3.37}
\sum_{i=2}^3\sum_{j=1}^3I_{ij}=\mathbb E\int_{Q_h} D_4|Y_h|^2{\rm d}t+\mathbb E\int_{Q_h^-} B_2|\mathbf D_h^+Y_h|^2{\rm d}t
+\mathbb{E} \int_{0}^{T}R_1{\rm d}t,
\end{align}
where
\begin{align*}
D_4=&s^{3}\lambda ^{3}\mathbf D_{h}(A_{1}A_{2})-s\lambda\mathbf D_h\left(A_{1}(-\partial_tl+s\lambda ^{2}A_{3}+s\lambda A_{4}-\Psi)\right)-\frac{1}{2} \Delta _{h}((1+A_0)\Psi)\\
&-s ^{2}\lambda ^{2}A_{2}\Psi +(-\partial_tl+s\lambda ^{2}A_{3}+s\lambda A_{4}-\Psi)\Psi,\\
B_2=&s\lambda\mathbf D _{h}^{+}(A_{1}(1+A_0))-\frac{h^{2}}{2}s^{3}\lambda ^{3}\mathbf D_{h}^{+}(A_{1}A_{2})+\frac{h^2}{2}s\lambda\mathbf D_h^+\left(A_{1}(-\partial_tl+s\lambda ^{2}A_{3}+s\lambda A_{4}-\Psi)\right)\\
&+\mathbf m_{h}^{+}\left((1+A_0)\Psi\right), \\
R_1=& s\lambda (A_{1}(1+A_0)(x_0,t))|(\mathbf D _{h}^{+}Y_{h})_{0}|^{2}-s\lambda (A_{1}(1+A_0)(x_{N+1},t))|(\mathbf D_{h}^{-}Y_{h})_{N+1}|^{2}.
\end{align*}

Therefore, from (\ref{3.36}) and (\ref{3.37}) it follows that
\begin{align}\label{3.38}
\sum_{i=1}^3\sum_{j=1}^3I_{ij}\geq&\frac{1}{4}\mathbb{E}\int_{Q_h^-}\left(\tau-2\right)s\lambda^2\varphi|\partial_x\psi|^2
\theta^{2}|g_h|^2{\rm d}t-C\mathbb E\int_{Q_h^-}s\varphi\theta^{2}|\mathbf D_{h}^{+}g_h|^{2}{\rm d}t\nonumber\\
&-\mathbb E\int_{Q_h}\theta^2 |f_{h}|^2{\rm d}t-C(\lambda)s^{2}e^{C(\lambda)s}||y_{h}(T)||_{L^2(\Omega,\mathcal F_T,\mathbb P;L^2(G_h))}^{2}\nonumber\\
&-\epsilon\mathbb E\int_{Q_h}\frac{1}{s\varphi}\theta^2 |\Delta _{h}y_{h}|^2{\rm d}t-C(\epsilon)\mathbb E\int_{Q_h}\theta^2 |\mathbf D_{h}^{+}y_{h}|^2{\rm d}t+\sum_{i=1}^4\mathbb{E} \int_{Q_h}D_i|Y_{h}|^{2}{\rm d}t\nonumber\\
&+\sum_{i=1}^2\mathbb E\int_{Q_h^-}B_i|\mathbf D_h^+ Y_h|^2{\rm d}t+\mathbb E\int_{Q_h^-}Ky_h\mathbf D_h^+y_h{\rm d}t+\mathbb{E} \int_{0}^{T}R_1{\rm d}t.
\end{align}

 \subsubsection{Positive lower bounds for the terms of $|Y_h|^2$ and $|\mathbf D_h^+ Y_h|^2$}
 
A direct calculation gives
\begin{align*}
D_1=&\frac{1}{2}s^{2}\lambda ^{2}(\partial_t f_2+\mathcal O_{\lambda}(sh))=s^2\mathcal O_{\lambda}(1),\\
D_2=&-\frac{1}{2} (-s\varphi_{tt}+s\lambda ^{2}(\partial_tf_3+\mathcal O_{\lambda}(sh))+s\lambda (\partial_tf_4+\mathcal O_{\lambda}(sh))-\tau s \varphi_{xxt} )=s\mathcal O_{\lambda}(1),\\
D_3=&-\frac{1}{2}\left(\partial_t(\mathbf D_h^+A_0)\theta \mathbf D_h^+\theta-\mathbf D_h^+A_0\partial_t\theta\mathbf D_h^+\theta+\mathbf D_h^+A_0\theta\partial_t(\mathbf D_h^+\theta)\right)\theta^{-2}=s^2\mathcal O_{\lambda}(sh),\\
D_4=&s^{3}\lambda ^{3}(\mathbf D_{h}A_{1}\mathbf m_hA_{2}+\mathbf m_h A_{1}\mathbf D_{h}A_{2}) -s ^{2}\lambda ^{2}A_{2}\Psi \\
&+(-\partial_tl+s\lambda ^{2}A_{3}+s\lambda A_{4}-\Psi)\Psi-s\lambda\mathbf D_h A_{1}\mathbf m_h(-\partial_tl+s\lambda ^{2}A_{3}+s\lambda A_{4}-\Psi)\\
&-s\lambda\mathbf m_h A_{1}\mathbf D_h(-\partial_tl+s\lambda ^{2}A_{3}+s\lambda A_{4}-\Psi)-\frac{1}{2} \Delta _{h}\Psi -\frac{1}{2}\Delta _{h}(A_0\Psi)\\
=&s^{3}\lambda ^{3}\left(\partial_x(\varphi\partial_x\psi)(\varphi^2|\partial_x\psi|^2)+\varphi\partial_x\psi\partial_x(\varphi^2|\partial_x\psi|^2)
+\mathcal O_{\lambda}(sh)\right)-\tau s ^{3}\lambda ^{2}\varphi^2\partial_{xx}\varphi|\partial_x\psi|^2\\
&+s^3\mathcal O_\lambda(sh)+s^2\mathcal O_\lambda(1)\nonumber\\
=&(3-\tau)\left(s^3\lambda^4\varphi^3|\partial_x\psi|^4+s^3\lambda^3\varphi^3\partial_{xx}\psi|\partial_x\psi|^2\right)+s^3\mathcal O_\lambda(sh)+s^2\mathcal O_\lambda(1)
\end{align*}
and
\begin{align*}
&\hspace{-0.8cm}B_1=h\mathcal O_{\lambda}(sh)+\mathcal O_{\lambda}(sh)=\mathcal O_{\lambda}(sh),\\
&\hspace{-0.8cm}B_2=\mathbf m_h^+\Psi+\mathbf m_h^+A_0\mathbf m_h^+\Psi+\frac{h^2}{4}\mathbf D_{h}^{+}A_0\mathbf D_h^+\Psi+s\lambda\left(\mathbf D _{h}^{+}A_{1}(1+\mathbf m_h^+A_0)+\mathbf m_h^+A_1\mathbf D_h^+A_0\right)
\end{align*}
\begin{align*}
&-\frac{h^{2}}{2}s^{3}\lambda ^{3}(\mathbf D_{h}^{+}A_{1}\mathbf m_h^+A_{2}+\mathbf m_h^+A_1\mathbf D_{h}^{+}A_{2})+\frac{h^2}{2}s\lambda\mathbf D_h^+\left(A_{1}(-\partial_tl+s\lambda ^{2}A_{3}+s\lambda A_{4}-\Psi)\right)\\
=&\tau  \left(s\lambda^2\varphi|\partial_x\psi|^2+s\lambda\varphi\psi_{xx}+s\mathcal O_{\lambda}(sh)\right)+s\mathcal O_{\lambda}(sh)+sh^2\mathcal \mathcal O_{\lambda}(sh)\\
&+(s\lambda^2\varphi|\partial_x\psi|^2+s\lambda\varphi\partial_{xx}+s\mathcal O_{\lambda}(sh))+s^3h^2\mathcal O_{\lambda}(1)\\
=&(\tau+1) \left(s\lambda^2\varphi|\partial_x\psi|^2+s\lambda\varphi\psi_{xx}\right)+s\mathcal O_{\lambda}(sh),
\end{align*}
where we have used
\begin{align*}
\left\{\begin{array}{l}
\partial_t(\mathbf D_h^+A_0)\theta \mathbf D_h^+\theta=\mathcal O_{\lambda}(sh)(s\lambda\varphi\partial_x\psi+s\mathcal O_{\lambda}(sh))\theta^2=s\mathcal O_{\lambda}(sh)\theta^2,\\
\mathbf D_h^+A_0\partial_t\theta\mathbf D_h^+\theta=\mathcal O_{\lambda}(sh)(s\lambda\varphi\partial_x\psi+s\mathcal O_{\lambda}(sh))s\lambda\psi_t\theta^2=s^2\mathcal O_{\lambda}(sh)\theta^2,\\
\mathbf D_h^+A_0\theta\partial_t(\mathbf D_h^+\theta)=\mathcal O_{\lambda}(sh)\left(s\lambda\partial_t\varphi\partial_x\psi+s^2\lambda^2\varphi^2\partial_t\psi\partial_x\psi+s^2\mathcal O_{\lambda}(sh)\right)\theta^2=s^2\mathcal O_{\lambda}(sh)\theta^2.
\end{array}
\right.
\end{align*}

Since $\tau\in (3/2,3)$ and $|\partial_x\psi|>0$ for all $x\in {\overline G}$, we can choose $\lambda$ sufficiently large to satisfy
\begin{align}\label{3.39}
&(3-\tau)\left(s^3\lambda^4\varphi^3|\partial_x\psi|^4+s^3\lambda^3\varphi^3\partial_{xx}\psi|\partial_x\psi|^2\right)+s^3\mathcal O_\lambda(sh)+s^2\mathcal O_\lambda(1)\nonumber\\
\geq & \frac{1}{2}(3-\tau)s^3\lambda^4\varphi^3|\partial_x\psi|^4\geq Cs^3\lambda^4\varphi^3
\end{align}
and
\begin{align}\label{3.40}
&(\tau+1) \left(s\lambda^2\varphi|\partial_x\psi|^2+s\lambda\varphi\psi_{xx}\right)+s\mathcal O_{\lambda}(sh)\geq \frac{\tau+1}{2}s\lambda^2\varphi|\partial_x\psi|^2\geq Cs\lambda^2\varphi.
\end{align}
Then, we have
\begin{align}\label{3.41}
&\sum_{i=1}^4\mathbb{E} \int_{Q_h}D_i|Y_{h}|^{2}{\rm d}t+\sum_{i=1}^2\mathbb E\int_{Q_h^-}B_i|\mathbf D_h^+ Y_h|^2{\rm d}t\nonumber\\
\geq&C\mathbb{E} \int_{Q_h}s^3\lambda^4\varphi^3|Y_{h}|^{2}{\rm d}t+C\int_{Q_h^-}s\lambda^2\varphi|\mathbf D^+Y_{h}|^{2}{\rm d}t.
\end{align}

\subsubsection{The remainder of the proof of Theorem 2.3}

Since $\partial_x\psi=2(x-x^*)>0$ for all $x\in  {\overline G}$, we obtain
\begin{align}\label{3.42}
R_1=& s\lambda(\varphi\partial_x\psi+\mathcal O_{\lambda}(sh))\left(|(\mathbf D _{h}^{+}Y_{h})_{0}|^{2}-|(\mathbf D_{h}^{-}Y_{h})_{N+1}|^2\right)\geq -Cs\lambda\varphi|(\mathbf D_{h}^{-}Y_{h})_{N+1}|^2,
\end{align}
which leads to
\begin{align}\label{3.43}
\mathbb{E} \int_{0}^{T}R_1{\rm d}t& \geq -C\mathbb{E} \int_{0}^{T}s\lambda\varphi|(\mathbf D_{h}^{-}Y_{h})_{N+1}|^2{\rm d}t.
\end{align}
Then, we substitute (\ref{3.41}) and (\ref{3.43}) into (\ref{3.38}) to yield 
\begin{align}\label{3.44}
&\mathbb{E} \int_{Q_h}s^3\lambda^4\varphi^3|Y_{h}|^{2}{\rm d}t+\int_{Q_h^-}s\lambda^2\varphi|\mathbf D^+Y_{h}|^{2}{\rm d}t+\mathbb{E}\int_{Q_h^-}s\lambda^2\varphi\theta^{2}|g_h|^2{\rm d}t\nonumber\\
\leq&C\mathbb E\int_{Q_h}II_1+C\mathbb E\int_{Q_h}\theta^2 |f_{h}|^2{\rm d}t+C\mathbb E\int_{Q_h^-}s\varphi\theta^{2}|\mathbf D_{h}^{+}g_h|^{2}{\rm d}t+C\mathbb{E} \int_{0}^{T}s\lambda\varphi|(\mathbf D_{h}^{-}Y_{h})_{N+1}|^2{\rm d}t\nonumber\\
&+C(\lambda ) s^2e^{C(\lambda)s}\|y_{h}(T)\|_{L^2(\Omega,\mathcal F_T,\mathbb P;L^2(G_h))}^{2}+\epsilon C\mathbb E\int_{Q_h}\frac{1}{s\varphi}\theta^2 |\Delta _{h}y_{h}|^2{\rm d}t\nonumber\\
&+C(\epsilon)\mathbb E\int_{Q_h}\theta^2 |\mathbf D_{h}^{+}y_{h}|^2{\rm d}t+\mathbb E\int_{Q_h^-}Ky_h\mathbf D_h^+y_h{\rm d}t.
\end{align}
On the other hand, by (\ref{3.7}) we have
\begin{align}\label{3.45}
\mathbb E\int_{Q_h}II_1=&\mathbb E\int_{Q_h}\theta I\left(f_h{\rm d}t+g_h{\rm d}B(t)\right)-\mathbb E\int_{Q_h}I^2{\rm d}t\nonumber\\
\leq &\frac{1}{2}\mathbb E\int_{Q_h}\theta^2|f_n|^2{\rm d}t-\frac{1}{2}\mathbb E\int_{Q_h}I^2{\rm d}t.
\end{align}
Moreover, we have
\begin{align}\label{3.46}
&\mathbf D_h^{+}Y_h=\mathbf m_{h}^{+}\theta\mathbf D_h^{+}y_h+\mathbf D_h^{+}\theta\mathbf m_{h}^{+} y_h\nonumber\\
=& (1+\mathcal O_\lambda(sh))\theta \mathbf D_{h}^{+}y_h+(s\lambda\varphi\psi_x+s\mathcal O_\lambda(sh))\theta\left(y_h+\frac{h}{2}\mathbf D_{h}^{+}y_h\right)\nonumber\\
=&\theta\mathbf D_h^{+}y_h+\mathcal O_\lambda(sh)\theta\mathbf D_h^{+}y_h+sO_\lambda(1)\theta y_h.
\end{align}
Therefore, by (\ref{3.44})-(\ref{3.46}) and choosing $\mathcal O_{\lambda}(sh)$ sufficiently small and $s$ sufficiently large to absorb the last two terms of $\mathbf D^+_h y_h$ and $y_h$ on the right-hand side of (\ref{3.46}), we find that
\begin{align}\label{3.47}
&\mathbb{E} \int_{Q_h}s^3\lambda^4\varphi^3\theta^2|y_{h}|^{2}{\rm d}t+\int_{Q_h^-}s\lambda^2\varphi\theta^2|\mathbf D^+y_{h}|^{2}{\rm d}t+\mathbb{E}\int_{Q_h^-}s\lambda^2\varphi|\theta^{2}|g_h|^2{\rm d}t+\mathbb E\int_{Q_h}I^2{\rm d}t\nonumber\\
\leq&C\mathbb E\int_{Q_h}\theta^2 |f_{h}|^2{\rm d}t+C\mathbb E\int_{Q_h^-}s\varphi\theta^{2}|\mathbf D_{h}^{+}g_h|^{2}{\rm d}t+C\mathbb{E} \int_{0}^{T}s\lambda\varphi\theta^2(x_{N+1},t)|(\mathbf D_{h}^{-}y_{h})_{N+1}|^2{\rm d}t\nonumber\\
&+C(\lambda ) s^2e^{C(\lambda)s}\|y_{h}(T)\|_{L^2(\Omega,\mathcal F_T,\mathbb P;L^2(G_h))}^{2}+\epsilon C\mathbb E\int_{Q_h}\frac{1}{s\varphi}\theta^2 |\Delta _{h}y_{h}|^2{\rm d}t\nonumber\\
&+C(\epsilon)\mathbb E\int_{Q_h}\theta^2 |\mathbf D_{h}^{+}y_{h}|^2{\rm d}t+\mathbb E\int_{Q_h^-}Ky_h\mathbf D_h^+y_h{\rm d}t.
\end{align}

The first-order difference term of $\mathbf D_{h}^{+}y_{h}$ on the right-hand side of (\ref{3.47}) can be absorbed by the first-order difference term on the left-hand side of (\ref{3.47}) if we choose $\lambda$ sufficiently large such that $\lambda\geq C(\epsilon)$.
To handle the second-order central difference term, we have to express  $\Delta_h Y_h$ in terms of $I$ and provide an estimate by the terms on the left-hand side of (\ref{3.47}).
The definition of $I$ gives
\begin{equation*}
	(1+A_{0})\Delta _{h}Y_{h}=-I+(-s^{2}\lambda^{2}A_{2}-\partial_t l+s\lambda^{2}A_{3}+s\lambda A_{4}-\Psi_{h} )Y_h.
\end{equation*}
Then, we obtain
\begin{align*}
		(1+A_{0})^2|\Delta _{h}Y_{h}|^2\leq& |I|^2+\left(-s^{2}\lambda^{2}A_{2}-\partial_t l+s\lambda^{2}A_{3}+s\lambda A_{4}-\Psi\right)^2 |Y_{h}|^2\nonumber\\
\leq &|I|^2+C\left(s^4\lambda^4\varphi^4|\partial_x\psi|^4+s^4\mathcal O_{\lambda}(sh)+s^2\mathcal O_{\lambda}(1)\right)|Y_{h}|^2,
\end{align*}
which implies 
\begin{align}\label{3.48}
\mathbb E\int_{Q_h}(1+A_{0})^2\frac{1}{s\varphi}|\Delta _{h}Y_{h}|^2{\rm d}t\leq \mathbb E\int_{Q_h}I^2{\rm d}t+C\mathbb E\int_{Q_h}s^3\lambda^4\varphi^3|Y_{h}|^2{\rm d}t
\end{align}
for  sufficiently large $s$  such that $s\varphi\geq 1$.
If we choose $h$ sufficiently small such that $|A_0|=|\mathcal O_{\lambda}(sh)|\leq \frac{1}{2}$, we have $(1+A_0)^2\geq \frac{1}{4}$. Then by (\ref{3.48}), we obtain 
\begin{align}\label{3.49}
\mathbb E\int_{Q_h}\frac{1}{s\varphi}|\Delta _{h}Y_{h}|^2{\rm d}t\leq C\mathbb E\int_{Q_h}I^2{\rm d}t+C\mathbb E\int_{Q_h}s^3\lambda^4\varphi^3|Y_{h}|^2{\rm d}t.
\end{align}
We use (\ref{2.1}), (\ref{2.3}) and (\ref{3.22}) to obtain
\begin{align}\label{3.50}
\Delta_hY_{h}=&\Delta_h\theta \mathbf m_{h}y_h+\mathbf m_{h}\theta\Delta_h y_h+2\mathbf D_h \theta\mathbf D_h y_h\nonumber\\
=&{\Delta_{h}\theta}y_{h}+\left(\theta+\frac{h^{2}}{2}\Delta_{h}\theta\right)\Delta_{h}y_{h}
+2{\mathbf D_{h}\theta}\left(\mathbf D_h^+y_{h}+\frac{h}{2}\Delta_h y_h\right)\nonumber\\
=&\left(s^2\lambda^2\varphi^2|\psi_x|^2+s\lambda\partial_x(\varphi\partial_x\psi)+s^2\mathcal O_{\lambda}(sh)\right)\theta y_h+\left(1+s^2h^2\mathcal O_{\lambda}(1)\right)\theta\Delta_h y_h\nonumber\\
&+2(s\lambda\varphi\psi+s\mathcal O_{\lambda}(sh))\theta\left(\mathbf D_h^+ y_{h}+\frac{h}{2}\Delta_h y_h\right)\nonumber\\
=&\theta\Delta_h y_h+\mathcal O_{\lambda}(sh)\theta\Delta_h y_h+s^2\lambda^2\varphi^2\mathcal O(1)\theta y_h+s\lambda\varphi\mathcal O(1)\theta\mathbf D_h^+ y_{h}.
\end{align}
Then, combining (\ref{3.49}) and (\ref{3.50}) we obtain 
\begin{align}\label{3.51}
&\epsilon\mathbb E\int_{Q_h}\frac{1}{s\varphi}\theta^2|\Delta _{h}y_{h}|^2{\rm d}t\nonumber\\
\leq& \epsilon C\mathbb E\int_{Q_h}\frac{1}{s\varphi}|\Delta _{h}Y_{h}|^2{\rm d}t+\epsilon C\mathbb{E} \int_{Q_h}s^3\lambda^4\varphi^3\theta^2|y_{h}|^{2}{\rm d}t+\epsilon C\int_{Q_h^-}s\lambda^2\varphi\theta^2|\mathbf D^+y_{h}|^{2}{\rm d}t\nonumber\\
\leq &\epsilon C\mathbb E\int_{Q_h}I^2{\rm d}t+\epsilon C\mathbb{E} \int_{Q_h}s^3\lambda^4\varphi^3\theta^2|y_{h}|^{2}{\rm d}t+\epsilon C\int_{Q_h^-}s\lambda^2\varphi\theta^2|\mathbf D^+y_{h}|^{2}{\rm d}t.
\end{align}
Next we deal with the last term on the right-hand side of (\ref{3.47}). By 
\begin{align*}
K=\mathbf D_h^+A_0\theta\partial_t\theta+h\mathbf D_h^+A_0\mathbf D_h^+\theta\partial_t\theta=s\mathcal O_{\lambda}(sh)\theta^2+s^2hO_{\lambda}(sh)\theta^2=s\mathcal O_{\lambda}(sh)\theta^2,
\end{align*}
 we have
\begin{align}\label{3.52}
&\mathbb E\int_{Q_h^-}Ky_h\mathbf D_h^+y_h{\rm d}t=\mathbb E\int_{Q_h^-}s\mathcal O_{\lambda}(sh)\theta^2y_h\mathbf D_h^+y_h{\rm d}t\nonumber\\
\leq &\mathbb E\int_{Q_h}s\mathcal O_{\lambda}(sh)\theta^2|y_h|^2{\rm d}t+\mathbb E\int_{Q_h^-}s\theta^2|\mathbf D_h^+y_h|^2{\rm d}t.
\end{align}

Therefore, by substituting (\ref{3.51}) and (\ref{3.52}) into (\ref{3.47}) we obtain
\begin{align}
&\mathbb{E} \int_{Q_h}\left((1-\epsilon C)s^3\lambda^4-s\mathcal O_{\lambda}(sh)\right)\varphi^3\theta^2|y_{h}|^{2}{\rm d}t+\int_{Q_h^-}\left((1-\epsilon C)s\lambda^2-s\right)\varphi\theta^2|\mathbf D^+y_{h}|^{2}{\rm d}t\nonumber\\
&+\mathbb{E}\int_{Q_h^-}s\lambda^2\varphi\theta^{2}|g_h|^2{\rm d}t+\int_{Q_h}(1-\epsilon C)I^2{\rm d}t\nonumber\\
\leq&C\mathbb E\int_{Q_h}\theta^2 |f_{h}|^2{\rm d}t+C\mathbb E\int_{Q_h^-}s\varphi\theta^{2}|\mathbf D_{h}^{+}g_h|^{2}{\rm d}t+C\mathbb{E} \int_{0}^{T}s\lambda\varphi\theta^2(x_{N+1},t)|(\mathbf D_{h}^{-}y_{h})_{N+1}|^2{\rm d}t\nonumber\\
&+C(\lambda ) s^2e^{C(\lambda)s}\|y_{h}(T)\|_{L^2(\Omega,\mathcal F_T,\mathbb P;L^2(G_h))}^{2}.
\end{align}
Finally, choosing $\epsilon$ sufficiently small, we can obtain the desired estimate (\ref{2.13}) and then complete the proof of Theorem 2.3.

\subsection{Proof of Theorem 2.4}
In order to apply Theorem 2.3 to prove (\ref{1-2.15}), we need to make the boundary conditions in (\ref{2.14}) homogeneous. To do this, we introduce $u_h$ satisfying\\
\begin{equation}
	\left\{\begin{array}{ll}
		{\rm d}u_h-\Delta _{h}u_{h}{\rm d}t=0,&(x_h,t)\in Q_{h}, \\
		(u_{h})_0=\gamma_1,\ (u_{h})_{N+1}=\gamma_2, &t\in(0,T),\\
          u_h(0)=0,&x\in G_h.
        \end{array}
	\right.
\end{equation}
Obviously, we know for a.e. $\omega\in \Omega$ that $u_h\in L^2(0,T;H^2(G_h))\cap C([0,T];L^2(G_h))$ satisfies
\begin{align}\label{3.55}
\|u_h\|_{L^2_{\mathcal F}(\Omega;L^2(0,T;H^2(G_h)))}+\|u_h\|_{L^2_{\mathcal F}(\Omega;C([0,T];L^2(G_h)))}\leq C\sum_{i=1}^2\|\gamma_i\|_{L^2_{\mathcal F}(\Omega;H^1(0,T))}.
\end{align}
Additionally, we can also obtain the following estimate for $\mathbf D^-_h u_h$ at boundary $x_{N+1}$:
\begin{align}\label{3.56}
\mathbb E\int_{0}^{T}|(\mathbf D_{h}^{-}u_{h})_{N+1}|^{2}{\rm d} t\leq C\sum_{i=1}^2\|\gamma_i\|^2_{L^2_{\mathcal F}(\Omega;H^1(0,T))}.
\end{align}
Indeed, we introduce a cut-off function $\chi\in C^\infty[0,L]$ such that $0\leq \chi(x)\leq 1$ and
\begin{equation}\label{3.57}
	\chi(x)=\left\{\begin{aligned}
		0,\quad&x\in[x_0, x_1],\\
		1, \quad&x\in[x_{N}, x_{N+1}].
	\end{aligned}\right.
\end{equation}
Letting $\hat u_h=\chi u_h$,  we easily see that 
\begin{equation}\label{3.58}
	\left\{\begin{array}{ll}
		{\rm d}\hat u_h-\Delta _{h}\hat u_{h}{\rm d}t=\hat f{\rm d}t,&(x_h,t)\in Q_{h}, \\
		(\hat u_{h})_0=0,\ (\hat u_{h})_{N+1}=\gamma_2, &t\in(0,T),\\
          \hat u_h(0)=0,&x\in G_h,
        \end{array}
	\right.
\end{equation}
where
\begin{align*}
\hat f=-\Delta_{h}\chi\mathbf m _{h}u_{h}-2\mathbf D_{h}\chi\mathbf D_{h}u_{h}-\frac{h^{2}}{4}\Delta _{h}\chi\Delta _{h}u_{h}
\end{align*}
satisfies
\begin{align}\label{3.59}
|\hat f|\leq C\left(|u_h|+|\mathbf D_h^+u_h|+h|\Delta_h u_h|\right).
\end{align}
We multiply the equation of $\Delta_h \hat u_h$ by $\mathbf D_h^- \hat u_h$ and integrate over $Q_h$. Then by using discrete integration by parts (\ref{2.5}) and (\ref{3.57}), we obtain 
\begin{align}\label{3.60}
\int_{Q_h} \partial_t \hat u_h \mathbf D_h^-\hat u_h {\rm d}t+\int_{Q_h^-}\Delta_h\hat u_h \mathbf D_h^+ \hat u_h{\rm d}t-\int_0^T
\left|(\mathbf D_h^- u_h)_{N+1}\right|^2{\rm d}t=\int_{Q_h}\hat f \mathbf D_h^- \hat u_h {\rm d}t,
\end{align}
which implies
\begin{align}\label{3.61}
\mathbb E\int_0^T
\left|(\mathbf D_h^- u_h)_{N+1}\right|^2{\rm d}t\leq C\left(\|\hat u_h\|_{L^2(0,T;H^2(G_h))}^2+\|\hat f\|_{L^2(0,T;L^2(G_h))}^2\right).
\end{align}
Therefore, from (\ref{3.55}), (\ref{3.59}) and (\ref{3.61}) we deduce (\ref{3.56}).

Letting $z_h=y_h-u_h$, we then obtain
\begin{equation}\label{}
	\left\{\begin{array}{ll}
		{\rm d}z_{h}-\Delta _{h}z_{h}{\rm d}t=f_{h}{\rm d}t+g_{h}{\rm d}B(t),&(x_h,t)\in Q_{h}, \\
		z_{h}=0, &(x_h,t)\in\Sigma_{h},\\
        z_{h}(0)=0,&x_h\in G_{h}.
        \end{array}
	\right.
\end{equation}
Applying (\ref{2.13}) to $z_h$, we have
\begin{align}\label{}
		&\mathbb{E} \int_{Q_h}\frac{1}{s\varphi}\theta^{2}|\Delta_{h}z_{h}|^{2}{\rm d}t+ \mathbb{E}\int_{Q_h^-}s\lambda^{2}\varphi\theta^{2}|\partial_{h}^{+}z_{h}|^{2}{\rm d}t+\mathbb{E} \int_{Q_h}s^{3}\lambda^{4}\varphi^{3}\theta^{2}|z_{h}|^{2}{\rm d}t\nonumber\\
&+\mathbb{E}\int_{Q_h^-}s\lambda^{2}\varphi\theta^{2}|g_{h}|^{2}dt\nonumber\\
\le&C\mathbb{E}\int_{Q_h}\theta^{2}|f_{h}|^{2}{\rm d}t+C\mathbb{E}\int_{Q_h^-}s\varphi\theta^{2}|\partial_{h}^{+}g_{h}|^{2}{\rm d}t+C\mathbb{E}
\int_{0}^{T}s\lambda\varphi\theta^{2}(x_{N+1},t)|(\mathbf D_{h}^{-}z_{h})_{N+1}|^{2}{\rm d} t\nonumber\\
&+ C(\lambda)s^{2}e^{C(\lambda)s}||z_{h}(T)||_{L^2(\Omega,\mathcal F_T,\mathbb P;L^2(G_h))}^{2},
\end{align}
which implies
\begin{align}\label{3.64}
		&\mathbb{E} \int_{Q_h}\frac{1}{s\varphi}\theta^{2}|\Delta_{h}y_{h}|^{2}{\rm d}t+ \mathbb{E}\int_{Q_h^-}s\lambda^{2}\varphi\theta^{2}|\partial_{h}^{+}y_{h}|^{2}{\rm d}t+\mathbb{E} \int_{Q_h}s^{3}\lambda^{4}\varphi^{3}\theta^{2}|y_{h}|^{2}{\rm d}t\nonumber\\
&+\mathbb{E}\int_{Q_h^-}s\lambda^{2}\varphi\theta^{2}|g_{h}|^{2}dt\nonumber\\
\le&C\mathbb{E}\int_{Q_h}\theta^{2}|f_{h}|^{2}{\rm d}t+C\mathbb{E}\int_{Q_h^-}s\varphi\theta^{2}|\partial_{h}^{+}g_{h}|^{2}{\rm d}t+C\mathbb{E}
\int_{0}^{T}s\lambda\theta^{2}(x_{N+1},t)|(\mathbf D_{h}^{-}y_{h})_{N+1}|^{2}{\rm d} t\nonumber\\
&+C(\lambda)s^{2}e^{C(\lambda)s}||y_{h}(T)||_{L^2(\Omega,\mathcal F_T,\mathbb P;L^2(G_h))}^{2}+C\mathbb{E}
\int_{0}^{T}s\lambda\varphi\theta^{2}(x_{N+1},t)|(\mathbf D_{h}^{-}u_{h})_{N+1}|^{2}{\rm d} t\nonumber\\
&+C(\lambda)s^{3}e^{C(\lambda)s}\left(||u_{h}||^2_{L^2_\mathcal F(\Omega;L^2(0,T;H^2(G_h)))}+||u_{h}(T)||_{L^2(\Omega,\mathcal F_T,\mathbb P;L^2(G_h))}^{2}\right).
\end{align}
Substituting (\ref{3.55}) and (\ref{3.56}) into (\ref{3.64}), we obtain the desired estimate (\ref{1-2.15}) and then complete the proof of Theorem 2.4.

\section{Proof of Theorem 2.5}
\vspace{2mm}

In this section, we will prove the stability of our discrete inverse random source problem, i.e. Theorem 2.5.

\vspace{2mm}

\noindent{\bf Proof of Theorem 2.5.} Letting $\tilde{y}_{h}=y_{h}^{(1)}-y_{h}^{(2)}$ and $\tilde{g}_{h}=g_{h}^{(1)}-g_{h}^{(2)}$, we obtain
\begin{align}\label{4.1}
	\left\{
\begin{array}{ll}
		{\rm d}\tilde y_{h}-\Delta _{h}\tilde y_{h}{\rm d}t=\left(a_h\tilde y_h+b_h\mathbf D_h^+ \tilde y_h\right){\rm d}t+\tilde g_h{\rm d}B(t),&(x_h,t)\in Q_{h}, \\
{\tilde y}_{h}=0, &(x_h,t)\in\Sigma_{h},\\
\tilde y_{h}(0)=0,	&x_h\in G_{h}.
\end{array}
	\right.
\end{align} 
Then applying Theorem 2.3 to $\tilde{y}_{h}$, we obtain
\begin{align}\label{4.2}
		& \mathbb{E}\int_{Q_h^-}s\lambda^{2}\varphi\theta^{2}|\partial_{h}^{+}\tilde y_{h}|^{2}{\rm d}t+\mathbb{E} \int_{Q_h}s^{3}\lambda^{4}\varphi^{3}\theta^{2}|\tilde y_{h}|^{2}{\rm d}t+\mathbb{E}\int_{Q_h^-}s\lambda^{2}\varphi\theta^{2}|\tilde g_{h}|^{2}dt\nonumber\\
\le&C\mathbb{E}\int_{Q_h}\theta^{2}\left(|\tilde y_h|^{2}+|\mathbf D^+_h \tilde y_h|^{2}\right){\rm d}t+C\mathbb{E}\int_{Q_h^-}s\varphi\theta^{2}|\partial_{h}^{+}\tilde  g_{h}|^{2}{\rm d}t\nonumber\\
&+C\mathbb{E}
\int_{0}^{T}s\lambda\varphi\theta^{2}(x_{N+1},t)|(\mathbf D_{h}^{-}\tilde  y_{h})_{N+1}|^{2}{\rm d} t+ C(\lambda)s^{2}e^{C(\lambda)s}\|\tilde y_{h}(T)\|_{L^2(\Omega,\mathcal F_T,\mathbb P;L^2(G_h))}^{2}.
\end{align}
By means of (\ref{2.15}), and choosing $\lambda$ sufficiently large to absorb the first and second terms on the right-hand side of (\ref{4.2}), we obtain
\begin{align}\label{4.3}
		& \mathbb{E}\int_{Q_h^-}s\lambda^{2}\varphi\theta^{2}|\partial_{h}^{+}\tilde y_{h}|^{2}{\rm d}t+\mathbb{E} \int_{Q_h}s^{3}\lambda^{4}\varphi^{3}\theta^{2}|\tilde y_{h}|^{2}{\rm d}t+\mathbb{E}\int_{Q_h^-}s\lambda^{2}\varphi\theta^{2}|\tilde g_{h}|^{2}dt\nonumber\\
\le&C\mathbb{E}
\int_{0}^{T}s\lambda\varphi\theta^{2}(x_{N+1},t)|(\mathbf D_{h}^{-}\tilde  y_{h})_{N+1}|^{2}{\rm d} t+ C(\lambda)s^{2}e^{C(\lambda)s}\|\tilde y_{h}(T)\|_{L^2(\Omega,\mathcal F_T,\mathbb P;L^2(G_h))}^{2}.
\end{align}
From (\ref{4.3}), we immediately deduce (\ref{2.16}) and then complete the proof of Theorem 2.5. \hfill$\Box$

\section{Proof of Theorem 2.6}

\setcounter{equation}{0}

In this section, we prove the conditional stability for the discrete Cauchy problem, i.e. Theorem 2.6. In the proof we borrow some ideas from [\ref{Ya2009}].

\vspace{2mm}

\noindent{\bf Proof of Theorem 2.6.}\  In order to estimate the solution of (\ref{1.4}) in $G_{0,h}\times (\epsilon, T-\epsilon)$ by Cauchy data at $x_{N+1}$, we need to choose a suitable weight function $\varphi$. Let $\tilde G=(0,L+\delta)$ with $\delta>0$, and let $\hat G\subset\tilde G\setminus\overline{G}$. It is easy to find a function $d\in C^{2}(\tilde G)$ such that
 \begin{equation}
 \left\{
 \begin{array}{ll}
 		d(x)>0,&x\in\tilde G,\\
 		d(x)=0,& x\in \partial\tilde G,\\ 
 		|\partial_x d(x)|>0,&x\in G\subset {\tilde G}\setminus \hat G.
 \end{array}
 	\right.
 \end{equation}
 A typical form of function $d$ is $d(x)=x((L+\delta)-x)$ with $\delta>L$.
Then, since $G_0\subset\subset \tilde G$, we can choose a sufficiently large $N>1$ such that
\begin{equation}
	\begin{aligned}
		G_0\subset \left \{ x\ |\ x\in \tilde G,\ d(x)>\frac{4}{N}\|d\|_{L^\infty(\tilde G)} \right \} \cap G. 
	\end{aligned}
\end{equation}
Moreover, we choose a positive number $\beta$ such that
\begin{equation}\label{5.3}
	\begin{aligned}
		\beta\epsilon^{2}< \|d\|_{L^\infty(\tilde G)}<2\beta \epsilon^{2}.
	\end{aligned}
\end{equation}
We arbitrarily fix $t_{0}\in [\sqrt{2}\epsilon,T-\sqrt{2}\epsilon ]$. Meanwhile, we denote 
\begin{align*}
\psi(x,t)=d(x)-\beta(t-t_0)^2,\quad \varphi(x,t)=e^{\lambda\psi(x,t)} 
\end{align*}
with fixed large parameter $\lambda>0$. Let \begin{align*}\mu_{k}=e^{\lambda\left(\frac{k}{N}\|d\|_{L^\infty(\tilde G)}-\frac{\beta\epsilon^{2}}{N}\right)}
\end{align*} 
and 
\begin{equation}
		Q^{(k)}=\left \{ (x,t) \ |\ x\in \overline{G},\ \varphi (x,t)>\mu_{k} \right \},\quad k=1,2,3,4.
\end{equation}
Then, we can verify that 
\begin{equation}\label{5.5}
	\begin{aligned}
		 {G_0}\times \left(t_{0}-\frac{\epsilon}{\sqrt{N}},t_{0}+\frac{\epsilon}{\sqrt{N}}\right)\subset Q^{(k)}\subset \overline{G}\times(t_{0}-\sqrt{2}\epsilon,t_{0}+\sqrt{2}\epsilon),\quad k=1,2,3,4.
	\end{aligned}
\end{equation}

Let $Q^{(k)}_h=Q^{(k)}\cap Q_h$. In order to apply Theorem 2.4, we introduce a cut-off function $\tilde\chi\in C^{\infty}(Q)$ such that $0\le\tilde\chi\le1$ and
\begin{equation}\label{5.6}
	\tilde\chi(x,t)=\left\{\begin{aligned}
		1,\quad \varphi(x,t)>\mu_{3},
		\\0,\quad \varphi(x,t)<\mu_{2}.
	\end{aligned}\right.
\end{equation}
Then letting $z_{h}=\tilde\chi y_{h}$, we obtain  
\begin{align}
		{\rm d}z_{h}-\Delta _{h}z_{h}{\rm d}t=(a_{h}z_{h}+b_{h}\mathbf D_{h}^{+}z_{h}+F_{h})dt+c_hz_{h}{\rm d}B(t),\quad
		(x_h,t)\in {Q_h},
\end{align} 
where \begin{align*}
F_{h}=&\partial_{t}\tilde\chi y_{h}-\Delta_{h}\tilde\chi\mathbf m _{h}y_{h}-2\mathbf D_{h}\tilde\chi\mathbf D_{h}y_{h}-\frac{h^{2}}{4}\Delta _{h}\tilde\chi\Delta _{h}y_{h}\\
&-b_h\mathbf D_{h}^{+}\tilde\chi\mathbf m _{h}^{+} y_{h}-\frac{h}{2}b_h\mathbf D_{h}^{+}\tilde\chi\mathbf D_{h}^{+}y_{h}.
\end{align*}
From the definition of $\tilde\chi$, together with $d(x_0)=0$, we have $(x_0,t)\in Q\setminus Q^{(1)}$ and then $\tilde\chi(x_0,t)=0$.  We further have  
\be\label{5.7}(z_h)_0=0,\quad t\in (0,T).\ee 
Since $t_0\in \left[\sqrt{2}\epsilon,T-\sqrt{2}\epsilon\right]$, we have 
$$\max\left\{\psi(x,0),\psi(x,T)\right\}  \leq d(x)-2\beta \epsilon^2 \leq 0,$$ which implies $(G \times \{t=0\})\cup(G \times \{t=T\}) \subset {\rm Supp} (\tilde\chi)$ and then 
\begin{align}\label{5.8}
z_h(0)=0,\quad{\rm and}\quad z_h(T)=0.
\end{align} 

Applying (\ref{1-2.15}) to $z_h$ and noting (\ref{5.7}) and (\ref{5.8}), we obtain
\begin{align}\label{5.9}
		&\mathbb{E}  \int_{Q_{h}}\frac{1}{s\varphi}\theta^{2}|\Delta_{h}z_{h}|^{2}{\rm d}t+ \mathbb{E}\int_{Q_{h}^{-}}s\lambda^{2}\varphi\theta^{2}|\mathbf D_{h}^{+}z_{h}|^{2}{\rm d}t+\mathbb{E} \int_{Q_{h}}s^{3}\lambda^{4}\varphi^{3}\theta^{2}|z_{h}|^{2}{\rm d}t\nonumber\\
		\le&C\mathbb{E} \int_{Q_{h}}\theta^{2}|a_{h}z_{h}+b_{h}\mathbf D_{h}^{+}z_{h}+F_{h}|^{2}{\rm d}t+C\mathbb{E}\int_{Q_{h}^{-}}s\varphi\theta^{2}|\mathbf D_h^{+}(c_hz_{h})|^{2}{\rm d}t\nonumber\\
		&+C\mathbb{E}
\int_{0}^{T}s\lambda\varphi\theta^{2}(x_{N+1},t)|\eta|^{2}{\rm d} t+C(\lambda)s^3e^{C(\lambda)s}\|\xi\|^2_{L^2(\Omega;H^1(0,T))}\nonumber\\
\le&C\mathbb{E} \int_{Q_{h}}\theta^{2}|F_{h}|^2{\rm d}t+C\mathbb{E}\int_{Q_{h}^{-}}s\varphi\theta^{2}\left(1+\frac{h}{2}\right)^2\left(|z_{h}|^2+|\mathbf D_h^{+}z_{h}|^{2}\right){\rm d}t\nonumber\\
		&+C(\lambda)s^3e^{C(\lambda)s}\left(\|\xi\|^2_{L^2(\Omega;H^1(0,T))}+\|\eta\|^2_{L^2(\Omega;L^2(0,T))}\right)
\end{align}
Firstly, we  notice that the second term on the right-hand side of (\ref{5.9}) can be absorbed by the second and third terms on the left-side of (\ref{5.9}). Secondly,
by (\ref{5.5}), we have $\overline {G_0}\times \left(t_{0}-\frac{\epsilon}{\sqrt{N}},t_{0}+\frac{\epsilon}{\sqrt{N}}\right)\subset Q^{(4)}$ and then
\begin{align}\label{5.11}
\overline {G_{0,h}}\times \left(t_{0}-\frac{\epsilon}{\sqrt{N}},t_{0}+\frac{\epsilon}{\sqrt{N}}\right)\subset Q_h^{(4)}\subset Q_h
\end{align}
for sufficiently small $h\in (0,h_1)$. Then, noting (\ref{5.6}) and (\ref{5.11}) we see that
\begin{align}\label{1-5.12}
&\mathbb{E}  \int_{Q_{h}}\frac{1}{s\varphi}\theta^{2}|\Delta_{h}z_{h}|^{2}{\rm d}t+ \mathbb{E}\int_{Q_{h}^{-}}s\lambda^{2}\varphi\theta^{2}|\mathbf D_{h}^{+}z_{h}|^{2}{\rm d}t+\mathbb{E} \int_{Q_{h}}s^{3}\lambda^{4}\varphi^{3}\theta^{2}|z_{h}|^{2}{\rm d}t\nonumber\\
\geq&\mathbb{E}  \int_{t_{0}-\frac{\epsilon}{\sqrt{N}}}^{t_{0}+\frac{\epsilon}{\sqrt{N}}}\left(\int_{G_{0,h}}\frac{1}{s\varphi}\theta^{2}
|\Delta_{h}y_{h}|^{2}+ \int_{G^-_{0,h}}s\lambda^{2}\varphi\theta^{2}|\mathbf D_{h}^{+}y_{h}|^{2}+\int_{G_{0,h}}s^{3}\lambda^{4}\varphi^{3}\theta^{2}|y_{h}|^{2}\right){\rm d}t\nonumber\\
\geq &C(\lambda)\frac{1}{s}e^{2\mu_4s}\mathbb{E}  \int_{t_{0}-\frac{\epsilon}{\sqrt{N}}}^{t_{0}+\frac{\epsilon}{\sqrt{N}}}\left(\int_{G_{0,h}}|\Delta_{h}y_{h}|^{2}+ \int_{G^-_{0,h}}|\mathbf D_{h}^{+}y_{h}|^{2}+\int_{G_{0,h}}|y_{h}|^{2}\right){\rm d}t.
\end{align}
On the other hand, we can verify that there exists a sufficiently small $h_2$ such that 
\begin{align}\label{5.12}
\partial_t\tilde\chi(x_h,t)=\mathbf D^+_h\tilde\chi(x_h,t)=\mathbf D_h\tilde\chi(x_h,t)=\Delta_h\tilde\chi(x_h,t)=0,\quad (x_h,t)\in Q_h'
\end{align}
for all $h\in (0,h_2)$, where 
\begin{align*}
Q_h'=\left\{(x_h,t)\in Q_h\ | \ \varphi(x_h,t)>\mu_3+\frac{1}{2}(\mu_4-\mu_3)\right\}.
\end{align*}
Indeed, we only need to prove 
\begin{align}\label{5.13}
\tilde\chi(x_h\pm h,t)\equiv 1, \quad (x_h,t)\in Q_h'.
\end{align}
 Since $\varphi(x_h,t)>\mu_3+\frac{1}{2}(\mu_4-\mu_3)$, we have $$\varphi(x_h\pm h,t)=\varphi(x_h,t)\left(1+\mathcal O_\lambda(h)\right)>\mu_3+\frac{1}{2}(\mu_4-\mu_3)+\mathcal O_\lambda(h).$$ Then we can choose $h_2$ sufficiently small such that $\varphi(x_h\pm h,t)>\mu_3$ for all $h\in(0,h_2)$. By the definition of $\tilde\chi$ we obtain (\ref{5.13}).  We apply (\ref{2.1}) and (\ref{5.12}) to yield
 \begin{align}\label{1-5.15}
 &\mathbb{E} \int_{Q_{h}}\theta^{2}|F_{h}|^2{\rm d}t\nonumber\\
 \leq& Ce^{2s\left(\mu_3+\frac{1}{2}(\mu_4-\mu_3)\right)}\mathbb E\int_{Q_h\setminus Q_h'}\left(|y_h|^2+|\mathbf m_h y_h|^2+|\mathbf m_h^+ y_h|^2+|\mathbf D_h y_h|^2\right){\rm d}t\nonumber\\
 &+Ce^{2s\left(\mu_3+\frac{1}{2}(\mu_4-\mu_3)\right)}\mathbb E\int_{Q_h\setminus Q_h'}\mathcal O(h)\left(|\mathbf D_h^+ y_h|^2+|\Delta_h y_h|^2\right){\rm d}t\nonumber\\
 \leq &Ce^{s\left(\mu_3+\mu_4\right)}\mathbb E\int_{Q_h\setminus Q_h'}\left(|y_h|^2+|\mathbf D_h^+ y_h|^2+|\Delta_h y_h|^2\right){\rm d}t.
 \end{align}
Therefore, from (\ref{5.9}), (\ref{1-5.12}) and (\ref{1-5.15}) it follows that
\begin{align*}
&C(\lambda)\frac{1}{s}e^{2\mu_4s}\mathbb{E}  \int_{t_{0}-\frac{\epsilon}{\sqrt{N}}}^{t_{0}+\frac{\epsilon}{\sqrt{N}}}\left(\int_{G_{0,h}}|\Delta_{h}y_{h}|^{2}+ \int_{G^-_{0,h}}|\mathbf D_{h}^{+}y_{h}|^{2}+\int_{G_{0,h}}|y_{h}|^{2}\right){\rm d}t\nonumber\\
\leq &Ce^{s\left(\mu_3+\mu_4\right)}\|y_h\|^2_{L^2_\mathcal F(0,T;H^2(G_h))}+C(\lambda)s^3e^{C(\lambda)s}\left(\|\xi\|^2_{L^2(\Omega;H^1(0,T))}+\|\eta\|^2_{L^2(\Omega;L^2(0,T))}\right)
\end{align*}
which implies
\begin{align*}
&\mathbb{E}  \int_{t_{0}-\frac{\epsilon}{\sqrt{N}}}^{t_{0}+\frac{\epsilon}{\sqrt{N}}}\left(\int_{G_{0,h}}|\Delta_{h}y_{h}|^{2}+ \int_{G^-_{0,h}}|\mathbf D_{h}^{+}y_{h}|^{2}+\int_{G_{0,h}}|y_{h}|^{2}\right){\rm d}t\nonumber\\
\leq &C(\lambda)se^{-s\left(\mu_4-\mu_3\right)}\|y_h\|^2_{L^2_\mathcal F(0,T;H^2(G_h))}+C(\lambda)s^4e^{C(\lambda)s}\left(\|\xi\|^2_{L^2(\Omega;H^1(0,T))}+\|\eta\|^2_{L^2(\Omega;L^2(0,T))}\right)\\
\leq& e^{-\frac{1}{2}s\left(\mu_4-\mu_3\right)}\|y_h\|^2_{L^2_\mathcal F(0,T;H^2(G_h))}+e^{2C(\lambda)s}\left(\|\xi\|^2_{L^2(\Omega;H^1(0,T))}+\|\eta\|^2_{L^2(\Omega;L^2(0,T))}\right)
\end{align*}
for all sufficiently large $s$. Noticing that $\mu_4-\mu_3>0$, for fixed $\lambda$ by the standard argument we obtain 
\begin{align}\label{5.15}
		&\|y_{h}\|_{L^2_\mathcal F\left(t_{0}-\frac{\epsilon}{\sqrt{N}},t_{0}+\frac{\epsilon}{\sqrt{N}};H^2(G_{0,h})\right)}\le CM^\kappa\left(\|\xi\|_{L^2(\Omega;H^1(0,T))}+\|\eta\|_{L^2(\Omega;L^2(0,T))}\right)^{1-\kappa}
	\end{align}
with 
$$\kappa=\frac{C(\lambda)}{C(\lambda)+\frac{1}{4}(\mu_4-\mu_3)}\in (0,1).$$
Finally, in (\ref{5.15}) taking $t_0=\sqrt{2}\epsilon+\frac{j\epsilon}{\sqrt{N}}$, $j=0,1,2,\cdots,m$ such that 
\begin{align*}
\sqrt{2}\epsilon+\frac{m\epsilon}{\sqrt{N}}\leq T-\sqrt{2}\epsilon\leq \sqrt{2}\epsilon+\frac{m\epsilon}{\sqrt{N}}
\end{align*}
and summing up over $j$, we obtain the desired estimate (\ref{2.19}) with  replacing $\epsilon$ by $\sqrt{2}\epsilon$ and then complete the proof of Theorem 2.6 due to $\epsilon$ is arbitrary. \hfill$\Box$

\vskip 0.5cm

{\bf Acknowledgement.}
This first author is supported by NSFC (No.12171248).

 \vspace*{0.5cm}
 
\newcounter{cankao}
\begin{list}
	{[\arabic{cankao}]}{\usecounter{cankao}\itemsep=0cm} \centerline{\bf
		References} \vspace*{0.5cm} \small

\item\label{Bao2010} Bao G, Chow S N, Li P, et al. Numerical solution of an inverse medium scattering problem with a stochastic source[J]. Inverse Problems, 2010, 26(7): 074014.

\item\label{Bao2012} Bao G, Xu X. An inverse random source problem in quantifying the elastic modulus of nanomaterials[J]. Inverse Problems, 2012, 29(1): 015006.
    
\item\label{Bar2003} Barbu V, R\u{a}\c{s}canu A, Tessitore G. Carleman estimates and controllability of linear stochastic heat equations[J]. Applied Mathematics and Optimization, 2003, 47: 97-120.
	
	\item\label{Bau2015} Baudouin L, Ervedoza S, Osses A. Stability of an inverse problem for the discrete wave equation and convergence results[J]. Journal de Mathématiques Pures et Appliquées, 2015, 103(6): 1475-1522.
	
	\item\label{Bau2013} Baudouin L, Ervedoza S. Convergence of an inverse problem for a 1-D discrete wave equation[J]. SIAM Journal on Control and Optimization, 2013, 51(1): 556-598.

	\item\label{Boy2020} Boyer F, Hernández-Santamaría V. Carleman estimates for time-discrete parabolic equations and applications to controllability[J]. ESAIM: Control, Optimisation and Calculus of Variations, 2020, 26: 12.
	
	\item\label{Boy2010} Boyer F, Hubert F, Le Rousseau J. Discrete Carleman estimates for elliptic operators and uniform controllability of semi-discretized parabolic equations[J]. Journal de mathématiques pures et appliquées, 2010, 93(3): 240-276.
	
	\item\label{Boy2011} Boyer F, Hubert F, Le Rousseau J. Uniform controllability properties for space/time-discretized parabolic equations[J]. Numerische Mathematik, 2011, 118(4): 601-661.

	\item\label{Boyer2010} Boyer F, Hubert F, Rousseau J L. Discrete Carleman estimates for elliptic operators in arb-itrary dimension and applications[J]. SIAM journal on control and optimization, 2010, 48(8): 5357-5397.
	
	\item\label{Boy2014} Boyer F, Le Rousseau J. Carleman estimates for semi-discrete parabolic operators and application to the controllability of semi-linear semi-discrete parabolic equations[J]. Annales de l'Institut Henri Poincare (C) Non Linear Analysis, 2014, 31(5): 1035-1078.
	
	\item\label{Cas2021} Casanova P G, Hernández-Santamaría V. Carleman estimates and controllability results for fully discrete approximations of 1D parabolic equations[J]. Advances in Computational Mathematics, 2021, 47: 1-71.
	
	\item\label{Cer2022} Cerpa E, Lecaros R, Nguyen T N T, et al. Carleman estimates and controllability for a semi-discrete fourth-order parabolic equation[J]. Journal de Mathématiques Pures et Appliquées, 2022, 164: 93-130.

	 \item\label{Dou2023} Dou F, Lü P. Stability and regularization for ill-posed Cauchy problem of a stochastic parabolic differential equation[J]. arxiv preprint arxiv: 2308.15741, 2023.

\item\label{Fernández2006} Fernández-Cara E, Guerrero S. Global Carleman inequalities for parabolic systems and applications to controllability[J]. SIAM journal on control and optimization, 2006, 45(4): 1395-1446.
    
\item\label{Gao2018} Gao P. Global Carleman estimates for the linear stochastic Kuramoto–Sivashinsky equations and their applications[J]. Journal of Mathematical Analysis and Applications, 2018, 464(1): 725-748.

\item\label{Hernández2021} Hernández-Santamaría V. Controllability of a simplified time-discrete stabilized Kuramoto-Sivashinsky system[J]. arXiv preprint arXiv:2103.12238, 2021.
	

 \item\label{Imanuvilov2001} Imanuvilov O Y, Yamamoto M. Carleman estimate for a parabolic equation in a Sobolev space of negative order and their applications[J]. Control of Nonlinear Distributed Parameter Systems, 2001, 113.

\item\label{Imanuvilov2005} Imanuvilov O Y, Yamamoto M. Carleman estimates for the non-stationary Lamé system and the application to an inverse problem[J]. ESAIM: Control, Optimisation and Calculus of Variations, 2005, 11(1): 1-56.

\item\label{Isakov2006}    Isakov V. Inverse problems for partial differential equations[M]. New York: Springer, 2006.
	
	\item\label{Kli1992} Klibanov M V. Inverse problems and Carleman estimates[J]. Inverse problems, 1992, 8(4): 575.

\item\label{Kli2004} Klibanov M V, Timonov A A. Carleman estimates for coefficient inverse problems and numerical applications[M]. Walter de Gruyter, 2004.
	
	\item\label{Lec2023} Lecaros R, Ortega J H, Pérez A, et al. Discrete Calderón problem with partial data[J]. Inverse Problems, 2023, 39(3): 035001.
	
	\item\label{Lec2021} Lecaros R, Ortega J H, Pérez A. Stability estimate for the semi-discrete linearized Benjamin-Bona-Mahony equation[J]. ESAIM: Control, Optimisation and Calculus of Variations, 2021, 27: 93.
	
	\item\label{Liu2019} Liu X, Yu Y. Carleman estimates of some stochastic degenerate parabolic equations and application[J]. SIAM Journal on Control and Optimization, 2019, 57(5): 3527-3552.

	\item\label{Lu2012} Lü Q. Carleman estimate for stochastic parabolic equations and inverse stochastic parabolic problems[J]. Inverse Problems, 2012, 28(4): 045008.
	
	\item\label{Lu2023} Lü Q, Zhang X. Inverse problems for stochastic partial differential equations: Some progresses and open problems[J]. Numerical Algebra, Control and Optimization. doi:10.3934/naco.2023 014.

	\item\label{Ngu2012} Nguyen T. Carleman estimates for semi-discrete parabolic operators with a discontinuous diffusion coefficient and application to controllability[J]. arxiv preprint arxiv:1211.2061, 2012.
	
	\item\label{Per2022} Pérez Contreras A A. On discrete carleman estimates: applications to controllability, stability and inverse problems[J]. 2022.
	
\item\label{Rousseau2010} Rousseau J L, Robbiano L. Carleman estimate for elliptic operators with coefficients with jumps at an interface in arbitrary dimension and application to the null controllability of linear parabolic equations[J]. Archive for rational mechanics and analysis, 2010, 195(3): 953-990.
    
\item\label{Tang2009}   Tang S, Zhang X. Null controllability for forward and backward stochastic parabolic equations[J]. SIAM Journal on Control and Optimization, 2009, 48(4): 2191-2216.

	\item\label{Wu2020} Wu B, Chen Q, Wang Z. Carleman estimates for a stochastic degenerate parabolic equation and applications to null controllability and an inverse random source problem[J]. Inverse Problems, 2020, 36(7): 075014.

	\item\label{Wu2022} Wu B, Liu J. On the stability of recovering two sources and initial status in a stochastic hyperbolic-parabolic system[J]. Inverse Problems, 2022, 38(2): 025010.
	
	\item\label{Ya2009} Yamamoto M. Carleman estimates for parabolic equations and applications[J]. Inverse problems, 2009, 25(12): 123013.
	
	\item\label{Yuan2017} Yuan G. Conditional stability in determination of initial data for stochastic parabolic equations[J]. Inverse Problems, 2017, 33(3): 035014.

\item\label{Yuan2021} Yuan G. Inverse problems for stochastic parabolic equations with additive noise[J]. Journal of Inverse and Ill-posed Problems, 2021, 29(1): 93-108.
	
	\item\label{Zhang2022} Zhang W, Zhao Z. Convergence analysis of a coefficient inverse problem for the semi-discrete damped wave equation[J]. Applicable Analysis, 2022, 101(4): 1430-1455.
	
	\item\label{Zhao2023} Zhao Z, Zhang W. Stability of a coefficient inverse problem for the discrete Schrödinger equation and a convergence result[J]. Journal of Mathematical Analysis and Applications, 2023, 518(1): 126665.

\end{list}
\end{document}